\documentclass[reqno,pdf]{siamart0516-mr}

\usepackage{lipsum}
\usepackage{amsfonts}
\usepackage{graphicx}
\usepackage{epstopdf}
\usepackage{algorithmic}
\ifpdf
  \DeclareGraphicsExtensions{.eps,.pdf,.png,.jpg}
\else
  \DeclareGraphicsExtensions{.eps}
\fi
\newcommand{\TheTitle}{On the full space--time discretization of the generalized Stokes equations: The Dirichlet case} 
\newcommand{\TheAuthors}{S. Eckstein and M. R\r{u}{\v{z}}i\v{c}ka}

\headers{Space--time discretization of the generalized Stokes equations}{\TheAuthors}

\title{{\TheTitle}}

\author{
  Sarah Eckstein\thanks{Institute of Applied Mathematics, Albert-Ludwigs-University Freiburg, Eckerstr. 1, D-$79104$ Freiburg,
    Germany.}
  \and
  Michael R\r{u}{\v{z}}i\v{c}ka\thanks{Institute of Applied Mathematics, Albert-Ludwigs-University Freiburg, Eckerstr. 1, D-$79104$ Freiburg,
    Germany (\email{rose@mathematik.uni-freiburg.de}).}
}

\usepackage{amsopn}

\ifpdf
\hypersetup{
  pdftitle={\TheTitle},
  pdfauthor={\TheAuthors}
}
\fi

\numberwithin{equation}{section}

\newtheorem{Theorem}[equation]{Theorem}
\newtheorem{Definition}[equation]{Definition}
     \newtheorem{Corollary}[equation]{Corollary}
     \newtheorem{Proposition}[equation]{Proposition}
     \newtheorem{Lemma}[equation]{Lemma}
  \newtheorem{Remark}[equation]{Remark}
\newtheorem{Assumption}[equation]{Assumption}

\def\enddemo{{\unskip\nobreak\hfill\penalty50\hbox{}\nobreak
  \hfil\parfillskip=0pt\finalhyphendemerits=0\hbox{\hskip5pt\vrule
   width5pt height5pt depth0pt\hskip1pt}\par\smallskip}}
\def\demo#1{\noindent{\sc Proof \rm#1:\enspace}\nobreak\qwdtrue}
\newif\ifqwd\qwdtrue

\def\Xint#1{\mathchoice
   {\XXint\displaystyle\textstyle{#1}}%
   {\XXint\textstyle\scriptstyle{#1}}%
   {\XXint\scriptstyle\scriptscriptstyle{#1}}%
   {\XXint\scriptscriptstyle\scriptscriptstyle{#1}}%
   \!\int}
\def\XXint#1#2#3{{\setbox0=\hbox{$#1{#2#3}{\int}$}
     \vcenter{\hbox{$#2#3$}}\kern-.5\wd0}}

\def\dashint{\Xint-}

\newcommand{\fe}[1]{\begin{bf} #1 \end{bf}}

 \def\Ruzicka{R\r{u}{\v{z}}i\v{c}ka\,}

\def\feu{{\fe{u}}}

\def\fev{{\fe{v}}}
\def\D{{\fe{D}}}
\def\F{{\fe{F}}}
\def\S{{\fe{S}}}
\def\f{{\fe{f}}}
\def\Du{{\fe{Du}}}
\def\Dv{{\fe{Dv}}}
\def\FDu{{\fe{F(Du)}}}

\newcommand{\SDu}{\fe{S(Du)}}

\newcommand{\feU}{\fe{U}}

\def\pihdivu{{\Pi_h^{\operatorname{div}}\feu }}
\def\pihdiv{{\Pi_h^{\operatorname{div}} }}

\begin{document}

\maketitle

\begin{abstract}
  In this work we treat the space-time discretization of the generalized Stokes equations in the case of 
  Dirichlet boundary conditions.
  We 
  prove error estimates in the case $p\in[\frac{2d}{d+2},\infty)$ that are
  independent of the degeneracy parameter $\delta\in[0,\delta_0]$.
  For $p\leq 2$, our convergence rate is optimal.
\end{abstract}

\begin{keywords}
  generalized Stokes equations, space-time discretization, error estimates
\end{keywords}

\begin{AMS}
 65M15, 
 65M60,  
   76A05, 
 35Q35. 
\end{AMS}


\section{Introduction}
The purpose of this paper is to establish an error analysis for the space-time discretization of the generalized Stokes system
\begin{equation}\label{kap1pstokes}
 \begin{aligned}
  \partial_t\fe{u}-\fe{\operatorname{div}}\SDu+\nabla q&=\fe{f}\qquad&&\text{in}\, I\times\Omega,\\
  \fe{\operatorname{div}\feu}&=0\qquad&&\text{in}\, I\times\Omega,\\
  \fe{u}(0)&=\fe{u}_0 &&\text{in}\,\Omega\\
  \feu&=\mathbf 0\qquad&&\text{at}\, I\times\partial \Omega,\\
 \end{aligned}
\end{equation}
for given external body force $\fe{f}=(f_1,...,f_d)$  and initial velocity $\feu_0$,
where ${\Omega\subset\mathbb{R}^d,}$ ${d\geq 2}$, is a bounded, polygonal domain and $I=(0,T),\,T>0$, is a bounded time interval.
The unknown functions are the velocity field $\feu=(u_1,...,u_d)$ and the pressure $q$.
The function $\fe{S}$ is the extra stress tensor, whose structure is given by characteristic properties of the examined fluid.
 Here, $\fe{S}$ depends on the symmetric part of the gradient of $\feu$, $\fe{Du}:=\frac{1}{2}(\partial_iu_j+\partial_ju_i)_{i,j=1,...,d}$.
The special case $\fe{S}=\fe{Id}$, i.e. $-\operatorname{div}\SDu=-\Delta\feu$, which leads to the Stokes equations. 
In this work, we will consider a more general situation. A typical example is given by  
 \begin{equation}\label{kap1phi}
  \SDu:=\varphi'(|\Du|)\frac{\Du}{|\Du|},
 \end{equation}
 where $\varphi'(t):=(\delta+t)^{p-2}t$ for some $p\in(1,\infty),\,\delta\in[0,\delta_0]$. Note that our results 
  carry over to the case 
  $\SDu:=\psi'(|\Du|)\frac{\Du}{|\Du|},$ where $\psi$ is an $N$-function that fulfills the equivalence $\psi'(t)\sim\varphi'(t)$.

The system \cref{kap1pstokes} is a simplification of the generalized Navier-Stokes equations   
For a broader discussion of these models we refer to \cite{MR} and \cite{MNRR}.

Our goal is to present a complete analysis for the space-time discretization of the generalized Stokes system \cref{kap1pstokes}.
Our main result will be the error estimate 
  \begin{equation}\label{kap1mainerrest}
      \|\feu-\feU\|_{L^{\infty}(I,L^{2}(\Omega))}+
         \|\FDu-\fe{F(DU)}\|_{L^{2}(I,L^{2}(\Omega))}\leq c\left(\Delta t+h^{\min\{1,\frac{2}{p}\}}\right),
      \end{equation}   
for $p\in[\frac{2d}{d+2},\infty)$ (see \Cref{geserr1}). 

Let us summarize some previous results. 
For the special case $p=2$, i.e.~the Navier-Stokes equations, Heywood and Rannacher established in a series of papers
\cite{HR1, HR2, HR3, HR4} complete existence and regularity results 
as well as an error analysis for the time-space discretization.
 Regarding \cref{kap1pstokes}, one of the main difficulties lies in the treatment of the stress tensor $\fe{S}$.
 Barrett and Liu \cite{BL1} introduced a quasi-norm technique to prove error estimates for the $p$-Laplacian.
   Later on, they also treated the parabolic $p$-Laplacian \cite{BL2} and $p$-fluids \cite{BL3}.
    In \cite{DER}, 
   Diening, Ebmeyer and \Ruzicka 
  adapted this technique for $N$-functions and proved optimal error
  estimates for parabolic systems with $p$-structure.

 Considering the treatment of the generalized Navier-Stokes equations, there are various results for the case $p\leq 2$ given 
 periodic boundary conditions.
 In \cite{PR} Prohl and \Ruzicka proved a first result for the space-time discretization of the generalized Navier-Stokes equations for some 
  $p\in(p_0,2]$.
 In a series of papers together with Diening \cite{DPR, DPR06}, the authors improved the sub-optimal results for the time discretization and increased 
  the range of admissible $p$'s.
   In \cite{BDR09}, Berselli, Diening and \Ruzicka proved optimal error estimates for the time discretization of the generalized Navier-Stokes equations
    in the case $p\in(\frac{3}{2},2]$ and, together with Belenki, the authors also proved error estimates for the finite element 
      approximation of the stationary generalized Stokes system, cf. \cite{BBDR}.
  In \cite{BDR14}, Berselli, Diening and \Ruzicka were finally able to prove the optimal estimate 
     $E_{h,\Delta t}\leq c(h+\Delta t)$ for $p\in(\frac{3}{2},2]$. 

%
   Previous results for the generalized Navier-Stokes equations
   discretize first in time and then in space.  Therefore, spatial
   regularity of the semi-discrete solution is needed.  This
   regularity can so far only be obtained in the setting of periodic
   boundary conditions.  In this work, we discretize first in space
   and then in time, as in \cite{HR1, HR2}. Therefore, we need time regularity of the
   semi-discrete solution, which we are able to prove even in the
   setting of Dirichlet boundary conditions. Moreover, our treatment
   includes for the first time also the case $p>2$.
   For $p\in[\frac{2d}{d+2},2]$, our error estimates are optimal.
   The results of this paper are based on the PhD thesis of S. Eckstein, cf. \cite{Sarah}.

   {This paper is organized as follows:} In \Cref{sectechtools}, we
   provide the necessary technical tools. We introduce $N$-functions
   and operators with $N$-potential. This provides the technical tools
   needed for handling the stress tensor.  Moreover, we look into the
   finite element approximation of divergence-free fields. We
   introduce suitable function spaces, discuss several interpolation
   results as well as the discrete inf-sup condition, which is
   necessary for the spatial approximation of the generalized Stokes
   system.  In \Cref{spatapp}, we briefly discuss existence and
   regularity results for \cref{kap1pstokes}. Then we introduce the
   corresponding spatial approximation. In \Cref{sectionreguh}, we
   show existence and regularity for the spatial approximation
   $\feu_h$ of $\feu$. Choosing a suitable approximation $\feu_0^h$
   for the initial value $\feu_0$, we are able to prove time
   regularity of $\feu_h$. Afterwards, we derive error estimates for
   the spatial error.  \Cref{The fully discrete solution} treats the
   fully discretized solution. We consider an implicit scheme and show
   existence and regularity as well as error estimates. We show that
   for the above-mentioned choice for the initial value, we can
   finally prove the error estimate \cref{kap1mainerrest} for
   $p\in[\frac{2d}{d+2},\infty)$ (see \Cref{geserr1}).
  
\section{Technical Tools}\label{sectechtools}
\subsection{Function Spaces}\label{labelsectionnfct}
Let $\Omega\subset\mathbb{R}^d$, $d\ge 2$, 
be an open, bounded domain. By $L^p(\Omega)$ and $W^{k,p}(\Omega)$,
$1\leq p\leq\infty,\,k\in\mathbb{N}$, we denote the classical Lebesgue
and Sobolev spaces, respectively.  An element of a $d$-dimensional
function space is distinguished from a scalar function by bold print,
i.e. $\fe{u}=(u_1,...,u_d)\in W^{k,p}(\Omega)$ means
${u_i\in W^{k,p}(\Omega)}$, $i=1,...,d$. We also use bold print to
indicate tensor-valued functions.  We define $W_0^{k,p}(\Omega)$ as
the closure of compactly supported functions $w \in C_0^\infty(\Omega)$ with respect to
$\|\cdot\|_{W^{k,p}(\Omega)}$.  By $L^p_0(\Omega)$ we define the
subspace of $L^p(\Omega)$ consisting of functions with vanishing mean
value
$\langle w\rangle_\Omega:=\frac{1}{|\Omega|}\int_\Omega w \,dx=0$ and
$W_{0,\operatorname{div}}^{1,p}(\Omega)$ is defined as
$W_{0,\operatorname{div}}^{1,p}(\Omega):=\{\fe{w}\in
W_{0}^{1,p}(\Omega)\big|\,\operatorname{div} \fe{w}=0\ \text{a.e. in}\
\Omega\}.$
The space $L^p_{\operatorname{div}}(\Omega)$ is defined as the closure
of $C_{0,\operatorname{div}}^\infty(\Omega)$
with respect to the $L^p$-norm.  For a Banach space $X$, we denote by
 $L^p(I,X),\,p\in[1,\infty],$ the classical Bochner spaces,
 cf. \cite{GaGr}.

By $C,\, c$ we denote generic constants, which may change from line to
line.  We say that two functions $f$ and $g$ are equivalent and use
the notation $f\sim g $, if there exist constants $c,\,C\geq 0$ such
that $cf\leq g\leq Cf$.  For normed vector spaces $X$ we denote the dual
space by $X^*$ and the duality product between $f\in X^*$ and $u\in X$
by $\langle f,u\rangle_{X^*,X}:=f(u)$ or simply by
$\langle f,u\rangle$, if there is no risk of confusion.  We will use
the notation
 \[(f,u):=\int\limits_\Omega fu\,dx,\]
 whenever the right-hand side is well-defined.

 The scalar product of two vectors $\fe{u},\,\fe{v}\in\mathbb{R}^d$ is
 denoted by $\fe{u}\cdot\fe{v}$.  For a tensor
 $\fe{A}\in\mathbb{R}^{d\times d}$ we denote its symmetric part by
 $\fe{A}^{\operatorname{sym}}:=\frac{1}{2}(\fe{A}+\fe{A}^\top)\in\mathbb{R}^{d\times
   d}_{\operatorname{sym}}:=\{\fe{A}\in\mathbb{R}^{d\times
   d}\big|\,\fe{A}=\fe{A}^\top\}$.
 For $\fe{A},\ \fe{B}\in\mathbb{R}^{d\times d}$ we denote by
 $\fe{A:B}$ the component-wise inner product and $|\fe{A}|$ denotes
 the Hilbert-Schmidt norm. 

 We will also use Orlicz and Sobolev-Orlicz spaces, cf. \cite{Krasno}.
 To this end, we use $N$-functions
 $\psi:\mathbb{R}^{\geq 0}\to\mathbb{R}^{\geq 0}$, as defined in
 \cite{DR1}.  We denote by $\psi^*$ its complementary function.
We say that $\psi$ fulfills the $\Delta_2$-condition, if there exists a constant $c>0$, such that for all $t\geq 0$, there holds $ \psi(2t)\leq c\psi(t).$
By $\Delta_2(\psi)$ we denote the smallest such constant.
 In the following we work solely with $N$-functions $\psi$, such that 
$\psi$ and $\psi^*$ satisfy the $\Delta_2$-condition. Under this condition we have
\[\psi^*(\psi'(t))\sim\psi(t).\]
We denote by $L^\psi(\Omega)$ and $W^{1,\psi}(\Omega)$ the classical Orlicz and Sobolev-Orlicz spaces,
i.e., $f\in L^\psi(\Omega)$ if the modular
$\rho(f):=\int_\Omega \psi(|f|)\,dx$ is finite and $f \in W^{1,\psi}(\Omega)$ if $f,\nabla f \in L^\psi(\Omega)$.
Note that the dual space 
$(L^\psi(\Omega))^*$ can
be identified with the space $L^{\psi^*}(\Omega)$.

\subsection{Basic properties of the extra stress tensor}\label{Basic properties of the extra stress tensor}
 In the whole paper we
assume that the extra stress tensor $\fe{S}$ has $N$-potential, which will be defined now. A
detailed discussion and full proofs can be found in \cite{DR1}.

\begin{Definition}[Operators with $N$-Potential]\label{pot1}
 Let $\psi$ be an $N$-function. We say that the operator $\fe{S}:\mathbb{R}^{d\times d}\to \mathbb{R}^{d\times d}_{\operatorname{sym}}$ possesses $N$-potential $\psi$, if $\fe{S(0)=0}$ and if for all $\fe{P}\in\mathbb{R}^{d\times d}\setminus\{\fe{0}\}$ there holds
\begin{equation}\label{pot2}
 \fe{S(P)}=\fe{S}_\psi(\fe{P}):=\frac{\psi'(|\fe{P}^{\operatorname{sym}}|)}{|\fe{P}^{\operatorname{sym}}|}\fe{P}^{\operatorname{sym}}.
\end{equation}
\end{Definition}
We want to concentrate on a special $N$-function with $(p, \delta)$-structure, which is for $t\geq 0$ given by
\begin{equation}\label{defphi}
 \varphi(t)=\int_0^t\varphi'(s)\,ds\qquad\text{ with} \ \varphi'(s) := (\delta + s)^{p-2}s .\end{equation}
The function $\varphi$ satisfies uniformly in $t$ the important equivalence
\[\varphi'(t) t\sim\varphi''(t),\]
since
$\min \{1, p - 1\} (\delta+t)^{p-2} \leq \varphi''(t)\leq \max \{1, p
- 1\}(\delta+t)^{p-2}$.
Moreover, $\varphi$ satisfies the $\Delta_2$-condition with
$\Delta_2(\varphi) \leq c 2^{\max \{2,p\}}$, hence independent of
$\delta$.  This implies that, uniformly in $t$, $\delta$, we have
\[\varphi'(t) t\sim\varphi(t).\]
The conjugate function $\varphi^*$ satisfies $\varphi^*(t) \sim (\delta^{p-1} +t)^{p'-2}t^2$.
 Also $\varphi^*$
satisfies the $\Delta_2$-condition with $\Delta_2(\varphi^*)\leq c2^{\max\{2,p'\}}.$
If $\varphi$ is given by \cref{defphi} the spaces $L^\varphi(\Omega)$
and $L^p(\Omega)$ coincide with uniform equivalence of the
corresponding norms. The constants only depend on $p$ and $\Omega$.

Throughout this paper, we are going to assume
\begin{Assumption}\label{stresstensor}
 The stress tensor $\fe{S}:\mathbb{R}^{d\times d}\to \mathbb{R}^{d\times d}_{\operatorname{sym}}$ possesses $N$-potential
 $\varphi$, where $\varphi$ is given by \cref{defphi}, with $p\in(1,\infty)$ and $\delta\in[0,1]$. 
\end{Assumption}
\begin{Remark}\label{delta0}
 Throughout this paper, the assumption $\delta\in[0,1]$ can be replaced by the assumption $\delta\in[0,\delta_0]$ for given $\delta_0>0$. 
 The estimates will then depend on $\delta_0$.
\end{Remark}
\begin{Remark}\label{}
  All results of this paper remain true, if we replace
  \Cref{stresstensor} by the assumption that $\fe{S}$ has $N$-potential $\psi$ for an $N$-function $\psi$ such that $\psi\sim\varphi$.
\end{Remark}

For an $N$-function $\psi$, we define the family of shifted $N$-functions $\{\psi_a\}_{a\geq 0}$
for $t\geq 0$ by $\psi_a(t):=\int\limits_0^t\psi_a'(s)\,ds,$ where 
\begin{equation}\label{sh3}
\psi_a'(t):=\psi'(a+t)\frac{t}{a+t}. 
\end{equation}
For the $N$-function defined in \cref{defphi} we have that $\varphi_a(t)\sim(\delta+a+t)^{p-2}t^2$ and also 
$(\varphi_a)^*(t)\sim((\delta+a)^{p-1}+t)^{p'-2}t^2$. The families $\{\varphi_a\}_{a\geq 0}$ and $\{(\varphi_a)^*\}_{a\geq 0}$ satisfy the $\Delta_2$-condition
uniformly in $a\geq 0$, with $\Delta_2(\varphi_a)\leq
c2^{\max\{2,p\}}$ and $\Delta_2((\varphi_a)^*)\leq c2^{\max\{2,p'\}}$,
respectively.

We need the following refined version of Young's inequality, cf. \cite{DR1}: 
\begin{Lemma}[Young's inequality]\label{young2}
 Let $\psi$ be an $N$-function with $\Delta_2(\psi)<\infty$ and $\Delta_2(\psi^*)<\infty$. Then, for every $\varepsilon>0$, there exists $c_\varepsilon>0$ only depending on 
 $\varepsilon,\ \Delta_2(\psi),$ and $\Delta_2(\psi^*)$ such that for all $s,\, t,\, a\geq 0$
\begin{equation}\label{young5}
  st\leq \varepsilon \psi_a(s)+c_\varepsilon\psi_a^*(t)
 \end{equation}
 and 
 \begin{equation}\label{young6}
 s\psi_a'(t)+t\psi_a'(s)\leq\varepsilon\psi_a(s)+c_\varepsilon\psi_a(t).
  \end{equation}
\end{Lemma}
Closely related to the extra stress tensor $\fe{S}$ is the function
$\fe{F}:\mathbb{R}^{d\times d}\to\mathbb{R}^{d\times d}_{\operatorname{sym}}$ defined through
\begin{equation}\label{pot14}
\fe{F(P)}:=(\delta+|\fe{P}^{\operatorname{sym}}|)^{\frac{p-2}{2}}\fe{P}^{\operatorname{sym}}.
\end{equation}
The connection between $\fe{S,\ F}$, and $\{\varphi_a\}_{a\geq 0}$ is best explained by the following
lemma (cf. \cite[Lemma 6.16]{DR1}).

\begin{Lemma}\label{pot15}
Let \Cref{stresstensor} be fulfilled and let $\fe{F}$ be defined as in \cref{pot14}. Then there holds
 for all $\fe{P,\,Q}\in\mathbb{R}^{d\times d}$
\begin{equation}\label{pot15a}
 \begin{split}
  \fe{(S(P)-S(Q)):(P-Q)}
    &\sim \fe{|F(P)-F(Q)|}^2\\
    &\sim \varphi_{|\fe{P}^{\operatorname{sym}}|}(|\fe{P}^{\operatorname{sym}}-\fe{Q}^{\operatorname{sym}}|)\\ 
 &\sim \varphi_{|\fe{Q}^{\operatorname{sym}}|}(|\fe{P}^{\operatorname{sym}}-\fe{Q}^{\operatorname{sym}}|)\\ 
 &\sim\varphi''(|\fe{P}^{\operatorname{sym}}|+|\fe{Q}^{\operatorname{sym}}|)|\fe{P}^{\operatorname{sym}}-\fe{Q}^{\operatorname{sym}}|^2,
 \end{split}
\end{equation}
where the constants only depend on $p$. 
Furthermore, we have
\begin{equation}\label{pot15b}
\fe{S(P):P\sim|F(P)|^2\sim}\varphi(|\fe{P}^{\operatorname{sym}}|). 
\end{equation}
\end{Lemma}
In this case we have 
\[
  (\SDu,\Du)\sim\int\limits_\Omega\varphi(|\Du|)\,dx\sim\int\limits_\Omega(\delta+|\Du|)^{p-2}|\Du|^2\,dx.
\]
 Moreover, the following estimate follows directly from \Cref{pot15} and Young's inequality \cref{young6}.
\begin{Lemma}\label{Au4}
 Let \Cref{stresstensor} be fulfilled and let $\fe{F}$ be defined as in \cref{pot14}.
 For all $\varepsilon>0$ exists $c_\varepsilon>0$, depending on $\varepsilon$ and the $\Delta_2$-constants such that for all sufficiently smooth vector fields $\fe{u,\, v,\, w},$ we have
 \begin{align*}
  \fe{(S(Du)-S(Dv),Dw-Dv)}
     &\leq\varepsilon\|\fe{F(Du)-F(Dv)}\|_2^2
       +c_\varepsilon\|\fe{F(Dw)-F(Dv)}\|_2^2\,,
\\
  \fe{(S(Du)-S(Dv),Dw-Dv)}
   &\leq\varepsilon\|\fe{F(Dw)-F(Dv)}\|_2^2 +c_\varepsilon\|\fe{F(Du)-F(Dv)}\|_2^2.
 \end{align*}
\end{Lemma}

\begin{Lemma}\label{Au6}
   Let $\Omega\subset\mathbb{R}^d$ be an open, bounded domain and let $\fe{S}$ fulfill \Cref{stresstensor}.
   For $p\in(1,2)$ we have
   \begin{equation}\label{Au7}
   \begin{split}
   \|\fe{F(Du)}&-\fe{F(Dv)}\|_2^{\frac{4}{p}} \leq\|\fe{Du-Dv}\|_p^2\\
       &\leq c\big(K+\|\fe{Du}\|_p+\|\fe{Du-Dv}\|_p\big)^{2-p}
   \|\fe{F(Du)-F(Dv)}\|_2^2
   \end{split}
   \end{equation}
    and for $p\in[2,\infty)$ we have
    \begin{equation}\label{Au8}
  \begin{split}
   \|\fe{Du}&-\fe{Dv}\|_p^p\leq\|\fe{F(Du)-F(Dv)}\|_2^2\\
     &\leq c\big(K+\|\fe{Du}\|_p+ \|\fe{Du-Dv}\|_p\big)^{p-2}\|\fe{Du-Dv}\|_p^2,
   \end{split}
   \end{equation}
   with constants $c$ independent of $\delta\in[0,1]$. The constant $K$ is given by
   \[K:=\delta |\Omega|^{\frac{1}{p}}\leq|\Omega|^{\frac{1}{p}}.\]
\end{Lemma}
\demo{} 
  See \cite[Lemma 2.2]{BL1} and \cite[Lemma 2.80]{Sarah}.\enddemo

\begin{Corollary}\label{Au9}
 Assume that  the assumptions of \Cref{Au6} are satisfied. Then ${\feu\in W^{1,p}(\Omega)}$ implies $\FDu\in L^2(\Omega)$. 
\end{Corollary}

We also use 

\begin{Theorem}\label{diffquot}
 Let $I\subset\mathbb{R}$ be an interval, $f\in L^p(I,L^q(\Omega)),\ 1<p,q<\infty$. Suppose there exists 
 a constant $K>0$ such that there holds $d_\tau f\in L^p(I',L^q(\Omega))$ and $\|d_\tau f\|_{L^p(I',L^q(\Omega))}\leq K$ for all 
 $0<\tau<\operatorname{dist}(I',\partial I)$. Then the weak derivative $\partial_t f$ exists and we have the 
 estimate $\|\partial_t f\|_{L^p(I,L^q(\Omega))}\leq K$.
\end{Theorem}
\demo{} The proof adapts the classical results for Sobolev spaces and can be found in \cite[Theorem 2.1]{Sarah}.\enddemo

\subsection{Finite Element Approximation}\label{Finite Element approximation}
For the spatial approximation of the generalized Stokes equations, we will need two different finite element spaces:
one for the approximation of the divergence-free velocity field and one for the approximation of the pressure. 
The choice of these two spaces is not arbitrary. In fact, they need to fulfill the {\it discrete inf-sup condition}. 

We start by explaining the triangulation of the domain.
 From now on we assume that $\Omega\subset\mathbb{R}^d$, $d\ge 2$, is a  bounded domain with polyhedral Lipschitz boundary. 
  Furthermore, we assume that for fixed $h>0$, $\mathcal{T}_h=\{T_i\}_{i=1,...,N}$ is a finite decomposition of $\Omega$ into simplices $T_i$. 
Let 
  $h_T:=\mathrm{diam}(T), h:=\max_{T\in\mathcal{T}_h} h_T$
 and let 
$\rho_T$
 denote the diameter of the largest closed ball contained in $\overline{T}$.
We assume that the mesh is such that any two elements of $\mathcal{T}_h$ meet only in entire common faces or sides or vortices (i.e. 
there are no hanging nodes), that the mesh is non-degenerate, i.e. there exists a constant $\sigma_0>0$ independent of $h$, such 
that $\max_{T\in\mathcal{T}_h}\frac{h_T}{\rho_T}\leq \sigma_0.$
Let $N_T$ be the neighborhood of $T$,
$N_T:=\bigcup\{ T'\big|\,\overline{T}\cap\overline{T'}\neq\emptyset\},$
and $S_T:=\mathrm{int}\big(\cup_{T'\in N_T}\overline{T'}\big)$. Under the above assumptions, it is clear that the number 
of simplices in every $S_T$ is bounded by a constant independent of $h_T$ and therefore
\begin{equation}\label{TST}
  |T|\sim h_T^d\sim|S_T|.
 \end{equation}
For $l\in\mathbb{N}_0$, let $\mathcal{P}_l(T)$ be the space of polynomials of degree less or equal to $l$ on $T$ and let
\[\mathcal{P}_l(\mathcal{T}_h):=\{v\in C^0(\overline{\Omega})\big|v\raisebox{-.5ex}{$|$}_{T}\in\mathcal{P}_l(T)\ \text{for all}\ T\in\mathcal{T}_h\}\]
be the space of piecewise polynomials.

Now we are ready to introduce suitable finite element spaces for the spatial approximation of the generalized Stokes equations. The natural setting for the continuous equation is to seek the velocity in $W_{0,\operatorname{div}}^{1,p}(\Omega)$ and the pressure in $L^{p'}_0(\Omega)$, hence our finite element spaces should be approximations of these spaces. We approximate $W_0^{1,p}(\Omega)$ by
\[X_h:=\{\fe{v}_h\in W_0^{1,p}(\Omega)\,\big|\,\fe{v}_h\raisebox{-.5ex}{$|$}_{T}\in\mathcal{P}_k(T)\ \forall\ T\in\mathcal{T}_h,\ \fe{v}_h
 \raisebox{-.5ex}{$|$}_{\partial\Omega}=0\}\]
and $L^{p'}(\Omega)$ by
 \[Y_h:=\{q_h\in L^{p'}(\Omega)\,\big|\,q_h\raisebox{-.5ex}{$|$}_{T}\in\mathcal{P}_r(T)\ \forall\ T\in\mathcal{T}_h\}\]
for some $k,\ r\in\mathbb{N}_0$. For the approximation of the pressure, we then define the space 
\[Q_h:= Y_h\cap L^{p'}_0(\Omega).\]
The discrete divergence-free space $V_h$ is then defined by
\[V_h:=\{\fe{v}_h\in X_h\big|\,( q_h,\fe{\operatorname{div}}\fe{v}_h)=0\ \forall\ q_h\in Y_h\}.\]
The choice of the polynomial degrees $k$ and $r$ is not arbitrary and plays an important role for the solvability of the discretized problem. The existence of stable pairings $X_h,\,Y_h$ has been widely discussed, see \cite{GR}.

The existence of interpolation operators for $X_h$ and $Y_h$ is quite
standard. Typical examples would be the Scott-Zhang operator \cite{SZ}
for $X_h$ and the Cl\'{e}ment operator \cite{BF} or also a version of
the Scott-Zhang operator for $Y_h$. However, we need to introduce
additional assumptions in order to guarantee well-posedness of the discretized problem as well as interpolation results in $V_h$ and 
 $Q_h$.

\begin{Assumption}\label{io1}
 Let $X_h,\, V_h,\,Y_h$ and $Q_h$ be defined as above with $\mathcal{P}_1(\mathcal{T}_h)\subset X_h$ and $\mathcal{P}_0(\mathcal{T}_h)\subset Y_h$.
We assume that there exist linear projection operators 
\begin{align*}
  \pihdiv&:\, W_0^{1,p}(\Omega)\to X_h,
\\
  \Pi_h^Y&:\, L^{p'}(\Omega)\to Y_h,
\end{align*}
   which fulfill the following assumptions
\begin{itemize}
 \item[(i)] $\pihdiv$ is divergence-preserving in the sense that for all $\fe{w}\in W_0^{1,p}(\Omega),\ \eta_h\in Y_h$
         \begin{equation}\label{io2}
         (\operatorname{div}\fe{w},\eta_h)=(\operatorname{div} \pihdiv\fe{w},\eta_h).
         \end{equation}

 \item[(ii)] $\pihdiv$ is locally $W^{1,1}$-stable in the sense that
   for all $\fe{w}\in W_0^{1,p}(\Omega),\ T\in\mathcal{T}_h$,
         \begin{equation}\label{io4}
          \dashint\limits_{T}|\pihdiv\fe{w}| \,dx\leq c\dashint\limits_{S_T}|\fe{w}| \,dx+c\dashint\limits_{S_T}h_T|\nabla\fe{w}| \,dx
           .
         \end{equation}
\item[(iii)]  There holds
      \begin{equation}\label{io4a}
        \pihdiv\fe{w}=\fe{w}\qquad\forall\ \fe{w}\in\mathcal{P}_{1}(\mathcal{T}_h).
      \end{equation}
 \item[(iv)] $\Pi_h^Y$ is locally $L^1$-stable in the sense that
   for all $ q\in L^{p'}(\Omega),\ T\in\mathcal{T}_h$,
         \begin{equation}\label{io5}
          \dashint\limits_T|\Pi_h^Y q| \,dx\leq c\dashint\limits_{S_T}|q| \,dx .
         \end{equation}
 \end{itemize}
\end{Assumption}
\begin{Remark}
 In order for \Cref{io1} to be fulfilled, there are only certain admissible pairings of polynomial
 degrees $k$ and $r$ for the spaces $X_h$ and $Y_h$, respectively, cf. \cite{GR}.
\end{Remark}

\begin{Remark}
 Note that since our choice of $X_h$ already includes zero boundary values, $\pihdiv$ also needs to preserve boundary values. The Scott--Zhang operator \cite{SZ} is one 
example for such an interpolation operator, but it needs to be modified in order to be divergence-preserving.  Following \cite{GL}, we show how this is done for the MINI element in three space dimensions.\\
 Let $\Omega\subset\mathbb{R}^{3}$.
 Assume that for every element $\eta_h\in Y_h$ holds $\eta_h\raisebox{-.5ex}{$|$}_{T}\in\mathcal{P}_1(T)$, $T\in\mathcal{T}_h$ and that for every $\fe{w}_h\in X_h,$ the restriction $\fe{w}_h\raisebox{-.5ex}{$|$}_{T}$ is the sum of a polynomial of $\mathcal{P}_1(T)$ and a bubble function 
  $b_T\in \mathcal{P}_4(T)\cap W_0^{1,p}(T)$.
 For each simplex $T\in\mathcal{T}_h$ we define the constant
  \[c_T:=\frac{\dashint_T\Pi_h^{SZ}\fe{w}-\fe{w} \,dx}{\dashint_T b_T \,dx},\] 
 where $\Pi_h^{SZ}$ is the Scott--Zhang operator. Now, we can show that the operator 
\[\pihdiv\fe{w}:=\Pi_h^{SZ}\fe{w}-\sum\limits_{T\in\mathcal{T}_h} c_Tb_T\]
  satisfies \Cref{io1}. Since $\Pi_h^{SZ}$ preserves zero boundary values and $b_T$ vanishes at every 
 edge of our triangulation and therefore particularly on $\partial\Omega$, $\pihdiv$ maps $W_0^{1,p}(\Omega)$ to $X_h$.
  The definition of $c_T$ gives the divergence-preserving property \cref{io2}: Since $\pihdiv\fe{w}-\fe{w}$ has zero 
  boundary values on $\Omega$ and $\nabla \eta_h\raisebox{-.5ex}{$|$}_{T}$ is constant for every $T\in\mathcal{T}_h$, we have
     \begin{equation*}
        \begin{split}
          \int\limits_\Omega\operatorname{div}(\pihdiv\fe{w}&-\fe{w})\eta_h \,dx 
        =-\int\limits_\Omega(\pihdiv\fe{w}-\fe{w})\cdot\nabla\eta_h \,dx\\
        &=-\sum\limits_{T\in\mathcal{T}_h}\int\limits_T(\pihdiv\fe{w}-\fe{w}) \,dx\cdot \nabla \eta_h\raisebox{-.5ex}{$|$}_{T}\\
          &=-\sum\limits_{T\in\mathcal{T}_h}\nabla \eta_h\raisebox{-.5ex}{$|$}_{T}\cdot\left(\int\limits_T(\Pi_h^{SZ}\fe{w}-\fe{w}) \,dx -\int\limits_T b_T \,dx
             \frac{\dashint_T\Pi_h^{SZ}\fe{w}-\fe{w} \,dx}{\dashint_T b_T \,dx}\right)\\
         &=0.
       \end{split}
        \end{equation*}
The local $W^{1,1}$-stability follows since the Scott--Zhang operator fulfills \cref{io4}. We have, using also \cref{TST},
  \begin{equation*}
   \begin{split}
    \dashint\limits_T|\pihdiv\fe{w}| \,dx&\leq\dashint\limits_T|\Pi_h^{SZ}\fe{w}| \,dx+|c_T|\big|\dashint\limits_T b_T \,dx\big|\\
     &\leq c\dashint\limits_{S_T}|\fe{w}| \,dx+c\dashint\limits_{S_T}h_T|\nabla\fe{w}| \,dx
       +\big|\dashint\limits_T\Pi_h^{SZ}\fe{w}-\fe{w} \,dx\big|\\
    &\leq c\dashint\limits_{S_T}|\fe{w}| \,dx+c\dashint\limits_{S_T}h_T|\nabla\fe{w}| \,dx,
   \end{split}
  \end{equation*}
 which is \cref{io4}. Other examples for $\pihdiv$ are given in \cite{GL}, \cite{BF}, \cite{BBDR} and \cite{GS}.  

 \end{Remark}

The results of \cite{DR2} stated in the theorem below clarify 
that \cref{io4} and \cref{io5} already provide sufficient approximability results.
\begin{Theorem}\label{io6}
Let $Z_h:=\{\fe{w}\in L^1_{loc}(\Omega)\,\big|\,\fe{w}\raisebox{-.5ex}{$|$}_{T}\in\mathcal{P}(T)\ \text{for all}\ T\in\mathcal{T}_h\}$ be a finite element space, where $\mathcal{P}_{r_0}(T)\subset\mathcal{P}(T)\subset\mathcal{P}_{r_1}(T)$ for $r_0\leq r_1\in\mathbb{N}_0$ and assume that for $l_0\in\mathbb{N}_0$ there exists an interpolation operator $\Pi_h:W^{l_0,1}(\Omega)\to Z_h$ such that
\begin{itemize}
 \item[a)] $\Pi_h$ is $W^{l,1}$-stable in the sense that for some $ l_0\leq l\leq r_0+1$ and $m\in\mathbb{N}_0$ holds uniformly in $T\in\mathcal{T}_h$ and $\fe{w}\in W^{l,1}(\Omega)$
 \begin{equation}\label{io7}
  \sum\limits_{j=0}^m\dashint\limits_T |h_T^j\nabla^j\Pi_h\fe{w}| \,dx\leq c(m,l)\sum\limits_{k=0}^lh_T^k\dashint\limits_{S_T}|\nabla^k\fe{w}| \,dx.
 \end{equation}
\item[b)] For all $\fe{w}\in\mathcal{P}_{r_0}(\mathcal{T}_h)$ holds $\Pi_h\fe{w}=\fe{w}$.
\end{itemize}
Let further $\Psi:[0,\infty)\to[0,\infty)$ be an $N$-function which satisfies the $\Delta_2$-condition.
Then there holds uniformly in $T\in\mathcal{T}_h$ and $\fe{w}\in W^{l,1}(\Omega)$
\begin{itemize}
 \item[(i)] $\Pi_h$ is Orlicz-stable in the sense that
       \begin{equation}\label{io8}
         \sum\limits_{j=0}^m\dashint_T\Psi(h^j_T|\nabla^j\Pi_h\fe{w}|) \,dx\leq c(m,l,\Delta_2(\Psi))\sum\limits_{k=0}^l\dashint_{S_T}\Psi(h^k_T|\nabla^k\fe{w}|) \,dx,
       \end{equation}
 \item[(ii)] $\Pi_h$ possesses an Orlicz-approximability property:
       \begin{equation}\label{io9}
         \sum\limits_{j=0}^l\dashint_T\Psi(h^j_T|\nabla^j(\fe{w}-\Pi_h\fe{w})|) \,dx\leq c(l,\Delta_2(\Psi),\sigma_0)\dashint_{S_T}\Psi(h^l_T|\nabla^l\fe{w}|) \,dx,
       \end{equation}
 \item[(iii)] $\Pi_h$ is Orlicz-continuous:
       \begin{equation}\label{io10}
         \dashint_T\Psi(h^l_T|\nabla^l\Pi_h\fe{w}|) \,dx\leq c(l,\Delta_2(\Psi),\sigma_0)\dashint_{S_T}\Psi(h^l_T|\nabla^l\fe{w}|) \,dx.
       \end{equation}
\end{itemize}
\end{Theorem}
\demo{}
The proof can be found in \cite{DR2}.
\enddemo

\begin{Remark}
\Cref{io1} guarantees that $\pihdiv$ fulfills the requirements of \Cref{io6} with $Z_h=X_h,\,r_0=1,\,m=0,\,l=1$ 
and that $\Pi_h^Y$ fulfills the requirements of \Cref{io6} with $r_0=l_0=l=m=0$. Moreover, due to the choice of $Y_h$, $\Pi_h^Y$ fulfills the requirements with $r_0=l_0=m=0$ and $l=1$.
Using inverse estimates, we can show that $\pihdiv$ also fulfills the requirements in the case $r_0=l_0=1,\,m=l=2$.  \end{Remark}

\noindent
 From \Cref{io6} we deduce
\begin{Lemma}\label{io11}
Let \Cref{io1} be fulfilled, let $\fe{S}$ satisfy \Cref{stresstensor} and let the associated operator $\fe{F}$ 
be defined as in \cref{pot14}. Then we have for $1<q<\infty$ and for
all sufficiently smooth enough vector fields $\fe{v}$
\begin{align}\label{io13}
  \|\fev-\pihdiv\fev\|_q +h\|\nabla\pihdiv\fe{v}\|_q 
  &\leq    ch\|\nabla\fev\|_q,
  \\
  \label{io14}\|\fev-\pihdiv\fev\|_q+h\|\nabla(\fev-\pihdiv\fev)\|_q &\leq 
                                                                       ch^2\|\nabla^2\fev\|_q,
  \\
  \label{io15}\|\fe{F(D}\pihdiv\fev)\|_2 &\leq c\|\fe{F(D\fev)}\|_2,
  \\
  \label{io16}\|\fe{F(Dv)-F(D}\pihdiv\fev)\|_2 &\leq ch\|\nabla\fe{F(Dv)}\|_2.
\end{align}
For the interpolation operator $\Pi_h^Y$ and for every $N$-function $\psi$ with $\Delta_2(\psi)<\infty$ we have
\begin{align}\label{io17}
  \int\limits_T\psi(|\Pi_h^Yq|) \,dx&\leq c\int\limits_{S_T}\psi(|q|) \,dx
  \\
  \label{io17a}
  \int\limits_T\psi(|q-\Pi_h^Yq|) \,dx&\leq c\int\limits_{S_T}\psi(h_T|\nabla q|) \,dx.
\end{align}
\end{Lemma}

\demo{} See \cite{DR2} and \cite{BBDR}.\enddemo
From \Cref{io1}  it can be shown that the existence of a divergence-preser\-ving interpolation operator as 
in \cref{io2} ensures that a discrete inf-sup condition is fulfilled by $X_h$ and $Q_h$. This is needed for the existence 
of a discrete pressure, cf. \Cref{discpress} below.

\begin{Lemma}\label{isfinal}
 Let \Cref{io1} be fulfilled and let $\varphi$ be defined by \cref{defphi} with $\delta\in[0,1],\,p\in(1,\infty)$. 
Then there exists a constant $c>0$, depending only on $p$ and
$\Omega$, 
such that 
\begin{equation}\label{is13}
   \|q_h\|_{L^{p'}(\Omega)}\leq c\underset{\boldsymbol{\xi}_h\in X_h\setminus\{0\}}{sup} \frac{(q_h,\operatorname{div}\boldsymbol{\xi}_h)}{\|\boldsymbol{\xi}_h\|_{W_0^{1,p}(\Omega)}}
\end{equation}
 holds for any $q_h\in Q_h$ and
\begin{equation}\label{is14}
 \int\limits_\Omega\varphi^*(|q_h|) \,dx\leq  \underset{\boldsymbol{\xi}_h\in X_h}{\sup}\Big( \int\limits_\Omega q_h\operatorname{div}\boldsymbol{\xi}_h \,dx-\frac{1}{c}\int\limits_\Omega \varphi(|\nabla\boldsymbol{\xi}_h|) \,dx\Big)
\end{equation}
 holds for any $q_h\in Q_h$.
 \end{Lemma}
\demo{} See \cite[Lemma 4.1]{BBDR}.\enddemo

   
   \section{Spatial Approximation}\label{spatapp}
\subsection{The continuous solution}
Before we discuss the spatial approximation of system \cref{kap1pstokes}, we will discuss 
some existence and regularity results for the continuous problem.
We assume that $\fe{S}$ fulfills \Cref{stresstensor}.

The first approach to show existence of a unique solution of \cref{kap1pstokes} is by using a Galerkin ansatz, solving the emerging ordinary differential equations,
establishing a priori estimates and then passing to the limit in the approximate system using monotone operator techniques and Minty's trick
\footnote{In view of Korn's inequality, this follows from standard monotone operator theory for evolution equations (cf. \cite{Ze2}) for 
$p\geq\frac{2d}{d+2}.$ The case $p>1$ can be handled as in \cite{DNR}. The case of generalized Navier-Stokes equations is treated in 
\cite{La} and \cite{Li}.}. For $p\in (1,\infty),\,\delta\in[0,1]$,
this leads to the existence of a unique weak solution $\feu$
satisfying the energy estimate 
\[\|\feu\|_{L^\infty(I,L^2(\Omega))}+\|\FDu\|_{L^2(I,L^2(\Omega))}\leq c\]
uniformly in $\delta\in[0,1]$, where the constant $c$ only depends on the data.

Using a completely different technique, namely linearization, maximal regularity, and a fixed point argument, Bothe und Pr\"{u}\ss\ 
\cite[Theorem 4.1]{BP} were able to show that for $\delta>0$ and $p\in(1,\infty)$, there exists a unique, strong 
solution on a maximal time interval, provided that the data is smooth
enough. More precisely, it is proved:
\begin{Theorem}\label{BPregularitytheorem}
 Let $\fe{f}\in L^r(J\times\Omega),\ \feu_0\in W^{2-\frac{2}{r}}(\Omega)\cap W_{0,\operatorname{div}}^{1,r}(\Omega)$ with $d+2<r<\infty$. 
 Moreover, assume that $\fe{S}$ fulfills \Cref{stresstensor} with $p\in(1,\infty)$ and $\delta\in(0,1]$. 
 Then, there exists a maximal time interval $I\subset J$ and a unique velocity field 
\[\feu\in L^r(I,W^{2,r}(\Omega))\cap W^{1,r}(I,L^r(\Omega))\]
 and a unique scalar function 
 \[q\in L^r(I,W^{1,r}(\Omega))\cap L^r(I,L^r_0(\Omega)), \]
 that solve \cref{kap1pstokes}.
\end{Theorem}
\begin{Remark}
 Under the same assumptions, Bothe and Pr\"{u}\ss\ \cite[Theorem 2.1]{BP} proved this result for the generalized Navier-Stokes equations. Their work also covers more general boundary conditions as well as other structures for $\varphi'$. However, the case $\delta=0 $ is not included.
\end{Remark}

\subsection{Existence and Regularity of the Finite Element Solution}\label{sectionreguh}

We now focus on the spatial approximation of the generalized Stokes system.
In \Cref{fe6}, we show existence of weak solutions $\feu_h$ of the spatial discretization. In order to estimate the error between space- and space-time 
approximation in \Cref{The fully discrete solution}, we need a certain
time regularity of $\feu_h$. This is accomplished for a special 
approximation $\feu_0^h$ for the initial value $\feu_0$ in \Cref{fe11}.

The weak formulation of \cref{kap1pstokes} suggests the discrete analogue: For a sufficiently smooth field 
$\fe{f}:I\times\Omega\to\mathbb{R}^d$ and $\feu_0^h\in V_h$ find $\feu_h\in C^1(\overline{I},X_h)$ such that for every $t\in \overline{I}$ there holds
   \begin{equation}\label{fe2}
   	\begin{split}
     (\partial_t\feu_h(t),\boldsymbol{\xi}_h )+(\fe{S(Du}_h(t)),\fe{D}\boldsymbol{\xi}_h)
     &=(\fe{f}(t),\boldsymbol{\xi}_h)\qquad\forall\,\boldsymbol{\xi}_h\in V_h,\\
       \feu_h(0)&=\feu_0^h\qquad\text{in}\,\Omega.
     \end{split}
   \end{equation}
 We will formulate the existence result and some a priori estimates in the next theorem.
\begin{Theorem}\label{fe6}
  Let $\fe{S}$ fulfill \Cref{stresstensor} with
  $\delta\in[0,1],\,p\in[\frac{2d}{d+2},\infty)$ and let \Cref{io1} be
  fulfilled. Suppose $\fe{f}\in W^{1,2}(I,L^2(\Omega))$ and
  $\feu_0\in W^{1,p}_{0,\operatorname{div}}(\Omega)$ and let   $M_1>0$ be such that
\begin{equation}\label{fe6ass1}
 \|\fe{f}\|_{W^{1,2}(I,L^2(\Omega))}+\|\feu_0\|_{W^{1,p}(\Omega)}\leq M_1.
\end{equation}
Let $\feu_0^h\in V_h$ be an approximation of $\feu_0$ such that there exists $M_2>0$ with
\begin{equation}\label{fe6ass2}
 \|\feu_0^h\|_{W^{1,p}(\Omega)}\leq M_2.
\end{equation}
Then there exists a unique solution $\feu_h\in C^1(\overline{I},V_h)$ of \cref{fe2}.
Furthermore, we have the estimates
\begin{align}\label{fe8}
\|\feu_h\|_{L^\infty(I,L^2(\Omega))}+\|\fe{F(Du}_h)\|_{L^2(I,L^2(\Omega))}&\leq c(M_1,M_2),\\
\label{fe9}
\|\partial_t\feu_h\|_{L^2(I,L^2(\Omega))}+\|\fe{F(Du}_h)\|_{L^\infty(I,L^2(\Omega))}&\leq c(M_1,M_2),
\end{align}
where the constants are independent of $h$.
\end{Theorem}
\demo{}
  The identity \cref{fe2} is a system of ordinary differential equation for $\feu_h$ which can be solved by standard methods.
  Let $\boldsymbol{\xi}_1,...,\boldsymbol{\xi}_N$ be a basis of $V_h$ and let $\alpha_1^0,\dots,\alpha_N^0\in\mathbb{R}$ be such that  
  \[\feu_0^h(x)=\sum\limits_{i=1}^N\alpha_i^0\boldsymbol{\xi}_i(x).\]
  Since the Gram matrix $G=((\boldsymbol{\xi}_i,\boldsymbol{\xi}_j))_{i,j=1,...,N}$ is invertible, P\'eano's theorem yields the existence of 
  a solution
  $\pmb{\alpha}^N(t):=(\alpha_1^N(t),\dots,\alpha_N^N(t))$ on an interval $[0,T^*],\ T^*\leq T$,
  of 
\begin{equation}\label{ode}
\begin{split}
\sum\limits_{i=1}^N\partial_t\alpha_i^N(t)(\boldsymbol{\xi}_i,\boldsymbol{\xi}_j)& =(\fe{f}(t),\boldsymbol{\xi}_j)-(\fe{S}(\sum\limits_{i=1}^N\alpha_i^N(t)\fe{D}\boldsymbol{\xi}_i),\fe{D}\boldsymbol{\xi}_j)\qquad\forall\,t\in [0,T^*],\\
\alpha_j(0)&=\alpha_j^0,
\end{split}
\end{equation}
for all $j=1,...,N$.
 This gives the solution of \cref{fe2} by defining 
  \[\feu_h(t,x):=\sum\limits_{i=1}^N\alpha_i^N(t)\boldsymbol{\xi}_i(x).\]  
 By choosing $\boldsymbol{\xi}_h=\feu_h(t)$, using \Cref{pot15} and the continuous Gronwall inequality, we get the a priori estimate
  \begin{equation}\label{odea}
\begin{split}
 \|\feu_h\|_{C(I^*,L^2(\Omega))}+\|\fe{F(Du}_h)\|_{L^2(I^*,L^2(\Omega))}
    &\leq c(\|\feu_0^h\|_{L^2(\Omega)},\,\|\fe{f}\|_{L^2(I,L^2(\Omega))})\\
    &\leq c(M_1,M_2),
\end{split}
\end{equation}
  where we used the embedding $W_0^{1,p}(\Omega)\hookrightarrow L^2(\Omega)$ for $p\geq\frac{2d}{d+2}$.
  Since the right-hand side of \cref{odea} is independent of $T^*$, we can extend the solution to the whole interval $I$, which also yields \cref{fe8}.
   The uniqueness of $\feu_h$ is easily shown using strict monotonicity.

   It remains to show the a priori estimate \cref{fe9}.
    Since $\fe{f}\in W^{1,2}(I,L^2(\Omega))$ we deduce
    that $\feu_h\in W^{2,2}(I,V_h)\hookrightarrow
    C^1(\overline{I},V_h)$.  Thus, 
    \[\partial_t\feu_h=\sum\limits_{i=1}^N\partial_t\alpha_i^N\boldsymbol{\xi}_i\in
    C^0(\overline I,V_h).\]
 We choose $\boldsymbol{\xi}_h=\partial_t\feu_h(s)$ in \cref{fe2} and integrate over $(0,t)$ to get
\[\int\limits_0^t\|\partial_t\feu_h(s)\|_2^2\,ds+\int\limits_0^t(\fe{S(Du}_h(s)),\fe{D}\partial_t\feu_h(s)) \,ds
\leq c\int\limits_0^T\|\fe{f}(s)\|_2^2 \,ds,\]
where we also used Young's inequality.
Due to our definition of $\fe{F}$ and $\fe{S}$
 we have
\[(\fe{S(Du}_h),\D\partial_t\feu_h)=(\fe{S(Du}_h),\partial_t\fe{D}\feu_h)\sim\frac{d}{dt}\|\fe{F(Du}_h)\|_{L^2(\Omega)}^2\]
and thus we deduce
\begin{equation*}
\begin{split}
    \int\limits_0^t\|\partial_t\feu_h(s)\|_{L^2(\Omega)}^2 \,ds+\|\fe{F(Du}_h(t))\|_{L^2(\Omega)}^2
      &\leq c\|\fe{f}\|_{L^2(I,L^2(\Omega))}^2+c\|\fe{F(Du}_0^h)\|_{L^2(\Omega)}^2\\
      &\leq c(M_1,M_2)
     \end{split}
\end{equation*}
for every $t\in(0,T)$. In the last step, we also used \Cref{Au9} to bound $\fe{F(Du}_0^h)$. This proves \cref{fe9}.
   \enddemo
\begin{Remark}\label{discpress}
 Once the existence of a solution $\feu_h\in C^1(\overline{I}, V_h)$ of \cref{fe2} is ensured, \Cref{isfinal} yields the existence of a discrete pressure $q_h\in C^0(\overline{I},Q_h)$ such that
\begin{equation*}
\begin{split}
(q_h(t),\operatorname{div}\boldsymbol{\xi}_h ) &= (\fe{f}(t),\boldsymbol{\xi}_h )
 + (\partial_t\feu_h(t),\boldsymbol{\xi}_h)+(\fe{S(Du}_h(t)),\fe{D}\boldsymbol{\xi}_h )\qquad\forall\,\boldsymbol{\xi}_h\in X_h.
 \end{split}\end{equation*}
 
\end{Remark}
In order to prove error estimates for the time discretization, we need higher time regularity for the solution $\feu_h$ of \cref{fe2}. To this end, we specify the initial value $\feu_0^h$. For given $\feu_0\in W_{0,\operatorname{div}}^{1,p}(\Omega)$, let $\feu_0^h\in V_h$ be the unique solution of
\begin{equation}\label{iv}
 (\fe{S(D}\feu_0^h),\fe{D}\boldsymbol{\xi}_h)=(\fe{S(D}\feu_0),\fe{D}\boldsymbol{\xi}_h)\qquad\text{for all}\ \boldsymbol{\xi}_h\in V_h.
\end{equation}
We have
\begin{Lemma}\label{iv1}
 For given $\feu_0\in W_{0,\operatorname{div}}^{1,p}(\Omega)$ with $\|\feu_0\|_{W^{1,p}(\Omega)}\leq M_1$, there exists a unique solution $\feu_0^h\in V_h$ of \cref{iv} such that
\begin{equation}\label{iv2}
 \|\fe{F(D}\feu_0^h)\|_{L^2(\Omega)}\leq c \|\fe{F(D}\feu_0)\|_{L^2(\Omega)}\leq c(M_1).
\end{equation}
Furthermore, for $p\in[\frac{2d}{d+2},\infty)$ we have
\begin{equation}\label{iv1a}
  \|\feu_0^h\|_{W^{1,p}(\Omega)}\leq c(M_1).
\end{equation}
\end{Lemma}
\demo{}
The existence of a function $\feu_0^h\in V_h$ that satisfies \cref{iv} follows from Brouwer's fixed point theorem.
By setting $\boldsymbol{\xi}_h=\feu_0^h$ in \cref{iv} and using \cref{pot15a} and \Cref{Au4} we get
\begin{align*}
 \begin{aligned}
  \|\fe{F(Du}_0^h)\|_{L^2(\Omega)}^2
  &\leq c(\fe{S(D}\feu_0^h),\fe{D}\feu_0^h)
  =c (\fe{S(D}\feu_0),\fe{D}\feu_0^h)\\
  &\leq \varepsilon \|\fe{F(Du}_0^h)\|_{L^2(\Omega)}^2 +c_\varepsilon \|\fe{F(Du}_0)\|_{L^2(\Omega)}^2,
 \end{aligned}
\end{align*}
which yields the first inequality of \cref{iv2} by choosing
$\varepsilon$ sufficiently small. The second part of \cref{iv2}
follows from \Cref{Au9}.

For $p\in[\frac{2d}{d+2},2)$ we use estimate \cref{Au7} with $\fev=\fe{0}$ and Young's inequality to get
\begin{align*}
 \begin{aligned}
  \|\fe{Du}_0^h\|_{L^p(\Omega)}^2 &\leq c\big(\|\fe{Du}_0^h\|_{L^p(\Omega)}+K\big)^{2-p}\|\fe{F(Du}_0^h)\|_{L^2(\Omega)}^2\\
  &\leq cK^{2-p}\|\fe{F(Du}_0^h)\|_{L^2(\Omega)}^2+c\|\fe{Du}_0^h\|_{L^p(\Omega)}^{2-p}\|\fe{F(Du}_0^h)\|_{L^2(\Omega)}^2\\
  &\leq cK^{2-p}\|\fe{F(Du}_0^h)\|_{L^2(\Omega)}^2 +\varepsilon \|\fe{Du}_0^h\|_{L^p(\Omega)}^2+ c_\varepsilon \|\fe{F(Du}_0^h)\|_{L^2(\Omega)}^{\frac{4}{p}}
 \end{aligned}
\end{align*}
where $K=\delta|\Omega|^{\frac{1}{p}}$. Choosing $\varepsilon$ sufficiently small we get together with \cref{iv2}
\begin{equation}\label{iv5}
 \|\fe{Du}_0^h\|_{L^p(\Omega)}^2\leq cK^{2-p}\|\fe{F(Du}_0)\|_{L^2(\Omega)}^2+ c \|\fe{F(Du}_0)\|_{L^2(\Omega)}^{\frac{4}{p}}.
\end{equation}
Now, Korn's inequality and \cref{iv2} yield \cref{iv1a}.

For the case $p\in[2,\infty)$, we use Korn's inequality and inequality \cref{Au8} to get
\[\|\fe{Du}_0^h\|_{L^p(\Omega)}^p\leq c \|\fe{F(Du}_0^h)\|_{L^2(\Omega)}^{2}\leq c \|\fe{F(Du}_0)\|_{L^2(\Omega)}^{2},\]
where we also used \cref{iv2}.
\enddemo

Using \cref{iv} as the definition for the initial value $\feu_0^h$, we can improve the regularity results for $\feu_h$.

\begin{Theorem}\label{fe11}
Let $\fe{S}$ fulfill \Cref{stresstensor} with
$\delta\in[0,1],\,p\in[\frac{2d}{d+2},\infty)$ and let \Cref{io1} be
fulfilled. Assume $\fe{f}\in W^{1,2}(I,L^2(\Omega))$ and $\feu_0\in
W^{1,p}_{0,\operatorname{div}}(\Omega)$ and let $M_1,\,M_3>0$ be such that
\begin{align}\label{fe11ass1}
 \|\fe{f}\|_{W^{1,2}(I,L^2(\Omega))}+\|\feu_0\|_{W^{1,p}(\Omega)}&\leq M_1
\\
\label{fe11ass2}
 \|\nabla\fe{S(D}\feu_0)\|_{L^2(\Omega)}&\leq M_3.
\end{align}
Also, let $\feu_0^h\in V_h$ be given as the solution of \cref{iv}.
Then there exists a unique solution $\feu_h\in C^1(\overline{I},V_h)$ of \cref{fe2}.
Furthermore, $\feu_h$ satisfies the estimates
\begin{align}\label{fe11a}
\|\feu_h\|_{L^\infty(I,L^2(\Omega))}+\|\fe{F(Du}_h)\|_{L^2(I,L^2(\Omega))}&\leq c(M_1),
\\
  \label{fe11b}
\|\partial_t\feu_h\|_{L^2(I,L^2(\Omega))}+\|\fe{F(Du}_h)\|_{L^\infty(I,L^2(\Omega))}&\leq c(M_1),
\\
\label{fe11c}
 \|\partial_t\feu_h\|_{L^\infty(I,L^2(\Omega))}+\|\partial_t\fe{F(Du}_h)\|_{L^2(I,L^2(\Omega))}&\leq c(M_1,M_3), 
 \end{align}
 where the constants are independent of the parameter $h$.
 \end{Theorem}
\demo{}
The existence of a solution as well as the a priori estimates \cref{fe11a} and \cref{fe11b} are shown in \Cref{fe6} and \Cref{iv1}. It remains to show \cref{fe11c}.
  To this end, we are using the othogonal projection $P_h: L^2(\Omega)\to V_h$ defined by $(P_h\fe{v},\boldsymbol{\xi}_h)=(\fe{v},\boldsymbol{\xi}_h)$ for all $\boldsymbol{\xi}_h\in V_h$, $\fe{v}\in L^2(\Omega)$. It is clear that $P_h$ is a self-adjoint, continuous projection and fulfills $\|P_h\fe{v}\|_2\leq\|\fe{v}\|_2$ for all $\fe{v}\in L^2(\Omega)$.
  At first, we prove that $\|\partial_t\feu_h(0)\|_{L^2(\Omega)}$ is uniformly bounded. 
  Since $\partial_t\feu_h(0)\in V_h$ it follows that $P_h(\partial_t\feu_h(0))=\partial_t\feu_h(0)$.
   Thus we get
\begin{equation}\label{fe11e}
 \begin{split}
  \|\partial_t\feu_h(0)\|_{L^2(\Omega)}&=\underset{\underset{\|\boldsymbol{\xi}\|_{L^2(\Omega)}\leq 1}{\boldsymbol{\xi}\in L^2(\Omega)}}{\sup} (\partial_t\feu_h(0),\boldsymbol{\xi} )
  =\underset{\underset{\|\boldsymbol{\xi}\|_{L^2(\Omega)}\leq 1}{\boldsymbol{\xi}\in L^2(\Omega)}}{\sup} (\partial_t\feu_h(0),P_h\boldsymbol{\xi} )\\
   &\leq\underset{\underset{\|\boldsymbol{\xi}_h\|_{L^2(\Omega)}\leq 1}{\boldsymbol{\xi}_h\in V_h}}{\sup} (\partial_t\feu_h(0),\boldsymbol{\xi}_h).\\
 \end{split}
\end{equation}
Next, we use \cref{fe2} at time $t=0$ and the definition of $\feu_0^h$, \cref{iv}, to get from \cref{fe11e}
\begin{equation}\label{troubleterm}
 \begin{split}
  \|\partial_t&\feu_h(0)\|_{L^2(\Omega)}
   \leq\underset{\underset{\|\boldsymbol{\xi}_h\|_{L^2(\Omega)}\leq 1}{\boldsymbol{\xi}_h\in V_h}}{\sup} 
    (\fe{f}(0),\boldsymbol{\xi}_h)-(\fe{S(Du}_0^h),\fe{D}\boldsymbol{\xi}_h)\\
   &\leq\underset{\underset{\|\boldsymbol{\xi}_h\|_{L^2(\Omega)}\leq 1}{\boldsymbol{\xi}_h\in V_h}}{\sup} 
     (\fe{f}(0),\boldsymbol{\xi}_h)-(\fe{S(D}\feu_0),\fe{D}\boldsymbol{\xi}_h)\\
   &\leq \underset{\underset{\|\boldsymbol{\xi}_h\|_{L^2(\Omega)}\leq 1}{\boldsymbol{\xi}_h\in V_h}}{\sup}
      \|\fe{f}\|_{L^\infty(I,L^2(\Omega))}\|\boldsymbol{\xi}_h\|_{L^2(\Omega)}+\|\nabla\fe{S(D}\feu_0)\|_{L^2(\Omega)}\|\boldsymbol{\xi}_h\|_{L^2(\Omega)}\\
   &\leq c(M_1,M_3). 
\end{split}
 \end{equation}
  Now we use difference quotients in order to prove \cref{fe11c}. For $t\in(0,T')$, where $T'<T$ and $0<\tau<T-T'$, 
 we take \cref{fe2} at time $t+\tau$, subtract \cref{fe2} at time $t$, divide by $\tau$ and choose the difference quotient
 $\boldsymbol{\xi}_h=d_\tau\feu_h(t)\in V_h$ as a test function in the emerging equation to get
\begin{align*}
 \begin{aligned}
  (&\partial_t d_\tau\feu_h(t), d_\tau\feu_h(t) )
  +\frac{1}{\tau^2}\big( \fe{S(Du}_h(t+\tau))-\fe{S(Du}_h(t)),\fe{Du}_h(t+\tau)-\fe{Du}_h(t) \big)\\
 &=( d_\tau\fe{f}(t), d_\tau\feu_h(t) ).
 \end{aligned}
\end{align*}
 \Cref{pot15}, Young's and H\"{o}lder's inequality yield
\[\frac{d}{dt}\|d_\tau\feu_h(t)\|_{L^2(\Omega)}^2+\|d_\tau\fe{F(Du}_h(t))\|_{L^2(\Omega)}^2\leq
    c\|d_\tau\fe{f}(t)\|_{L^2(\Omega)}^2+\|d_\tau\fe{u}_h(t)\|_{L^2(\Omega)}^2.\]
 Hence, the continuous Gronwall inequality, the boundedness of $\partial_t\feu_h(0)$ and \Cref{diffquot} imply the estimate
\begin{equation*}
\|\partial_t\feu_h\|_{L^\infty(I,L^2(\Omega))}^2+\|\partial_t\fe{F(Du}_h)\|_{L^2(I,L^2(\Omega))}^2\leq c\|\partial_t\fe{f}\|_{L^2(I,L^2(\Omega))}^2
  +c\|\partial_t\feu_h(0)\|_{L^2(\Omega)}^2.
\end{equation*}
This finally proves the a priori estimate \cref{fe11c}.
\enddemo
\subsection{Error Estimates for the Spatial Error}\label{sectionerrest}
The goal of this section is to finally prove error estimates between
the continuous solution $\feu$ of \cref{kap1pstokes} and its finite
element approximation $\feu_h$.
Motivated by the regularity results of \cite{BP} and \cite{DR3}, \cite{BDR} (in the space-periodic setting) we make the following assumption.

\begin{Assumption}\label{BP}
Let $\fe{S}$ fulfill \Cref{stresstensor} with $ \delta\in[0,1],\ 
 p\in[\frac{2d}{d+2},\infty).$
 Assume $\fe{f}\in W^{1,2}(I,L^2(\Omega))$, $\feu_0\in W^{1,p}_{0,\operatorname{div}}(\Omega),$ 
 and suppose
\begin{equation*}
\|\fe{f}\|_{W^{1,2}(I,L^2(\Omega))}+\|\feu_0\|_{W^{1,p}_{0}(\Omega)}+\|\nabla\fe{F(D}\feu_0)\|_{L^2(\Omega)}+\|\nabla\fe{S(D}\feu_0)\|_{L^2(\Omega)}\leq K_1. 
\end{equation*}
 Let the weak solution $(\feu,q)\in L^p(I,W_{0,\operatorname{div}}^{1,p}(\Omega))\times L^{p'}(I,L^{p'}_0(\Omega))$ of \cref{kap1pstokes} be such that
\begin{align*}
\|\feu\|_{W^{1,2}(I,L^2(\Omega))}+\|\feu\|_{L^2(I,W^{2,2}(\Omega))}+\|\feu\|_{L^p(I,W^{1,p}(\Omega))}+\|\F(\Du)\|_{W^{1,2}(I\times \Omega)}\leq K_2,\\
\| q\|_{L^{p'}(I,W^{1,p'}(\Omega))}\leq K_3.
\end{align*}
\end{Assumption}

\begin{Remark}\label{rem1}
 If we additionally assume $\f\in L^r(I\times\Omega)$, $\feu_0\in W^{r,2-\frac{2}{r}}(\Omega)$ for some $r>\max\{p,p',d+2\}$, and $\delta\in(0,1]$,
 then in \cite{BP} it is shown there exists a solution $(\feu,\,q)$ of \cref{kap1pstokes} with the property that 
 \[\|\feu\|_{L^r(I,W^{2,r}(\Omega))}+\|\feu\|_{W^{1,r}(I,L^r(\Omega))}+ \| q\|_{L^r(I,W^{1,r}(\Omega))}\leq c,\]
 see \Cref{BPregularitytheorem}. Similarly to \Cref{Au9} we can show that $\feu\in L^p(I,W^{1,p}(\Omega))$ implies $\FDu\in L^{2}(I,L^{2}(\Omega))$.
 For $p\geq 2$, we can show that $\nabla\F(\Du)$ is bounded in $L^2(I\times\Omega)$, which means that in this case, the existence of a solution $(\feu,\,q)$ fulfilling \Cref{BP} is ensured.
  For $p<2$, we have the estimate $\|\nabla\FDu\|_{L^{2}(I,L^{2}(\Omega))}\leq \delta^{\frac{p-2}{2}}\|\feu\|_{L^2(I,W^{2,2}(\Omega))}$. Together with the maximal regularity result from \cite{BP}, this also provides a solution $(\feu,q)$ fulfilling \Cref{BP} with a constant $K_2(\delta)$ that tends to infinity for $\delta\to 0$.
\end{Remark}
\begin{Remark}\label{rem2}
  Let \Cref{BP} be fulfilled and let $\feu_0^h\in V_h$ be given by \cref{iv}.
  Recall that in \Cref{fe11} we have shown that for the finite element solution $\feu_h$ the following estimates hold true
  \begin{equation}\label{rem2a}
  \begin{split}
  \|\feu_h\|_{L^\infty(I,L^2(\Omega))}+\|\fe{F(Du}_h)\|_{L^2(I,L^2(\Omega))}&\leq c(K_1),\\
  \|\partial_t\feu_h\|_{L^2(I,L^2(\Omega))}+\|\fe{F(Du}_h)\|_{L^\infty(I,L^2(\Omega))}&\leq c(K_1),\\
  \|\partial_t\feu_h\|_{L^\infty(I,L^2(\Omega))}+\|\partial_t\fe{F(Du}_h)\|_{L^2(I,L^2(\Omega))}&\leq c(K_1). 
  \end{split}
  \end{equation}
\end{Remark}

Let us now start estimating the error between $\feu$ and $\feu_h$. 
 As a first step, we prove a best approximation result.
\begin{Proposition}\label{bestapp}
  Let \Cref{BP} and \Cref{io1} be fulfilled.
   Moreover, assume that $\feu_0^h\in V_h$ is given by \cref{iv} and
   let $\feu_h$ be the corresponding finite element solution ensured
   by \Cref{fe11}.  
   Then we have\footnote{For the sake of readability, we omit the dependence on $t$ in \cref{bestapp1}.}
  \begin{equation}\label{bestapp1}
  \begin{split}
 \frac{d}{dt}\|&\feu-\feu_h\|_{L^2(\Omega)}^2+\|\fe{F(Du)}-\F(\Du_h)\|_{L^2(\Omega)}^2\\
     &\leq c\underset{\boldsymbol{\zeta}_h\in V_h}{\inf}\left(\|\partial_t(\feu-\feu_h)\|_{L^2(\Omega)}\|\feu-\boldsymbol{\zeta}_h\|_{L^2(\Omega)}
     +\|\fe{F(Du)}-\F(\D\boldsymbol{\zeta}_h)\|_{L^2(\Omega)}^2\right)\\ 
      &\ \ +c\underset{\mu_h\in Y_h}{\inf}\int\limits_\Omega\big(\varphi_{|\Du|}\big)^*(|q-\mu_h|) \,dx
  \end{split}
  \end{equation}
 for almost every $t\in I$.
\end{Proposition}
\demo{} 
  By subtracting \cref{fe2} from the weak formulation of \cref{kap1pstokes} and choosing $\boldsymbol{\xi}_h=\boldsymbol{\zeta}_h-\feu_h(t)\in V_h$ for arbitrary $\boldsymbol{\zeta}_h\in V_h$  as a test function, we get the error equation
  \begin{equation*}\begin{split}
  (&\partial_t(\feu-\feu_h)(t),\boldsymbol{\zeta}_h-\feu_h) +(\fe{S(Du(}t))-\fe{S(Du}_h(t)),\D(\boldsymbol{\zeta}_h-\feu_h)) \\
   &=(q(t),\operatorname{div}(\boldsymbol{\zeta}_h-\feu_h)) 
   =(q(t)-\mu_h,\operatorname{div}(\boldsymbol{\zeta}_h-\feu_h)) 
  \end{split}
  \end{equation*}
  for all $\boldsymbol{\zeta}_h\in V_h,\ \mu_h\in Y_h$ and almost every $t\in
  I$. In the last step we used that due to the definition of $V_h$ we
  have $(\mu_h,\operatorname{div}\boldsymbol{\xi}_h )=0$ for all $\boldsymbol{\xi}_h\in V_h,\ \mu_h\in Y_h$. 
     After rearranging the terms and using \Cref{pot15}, we obtain
   \begin{equation}\label{bestapp2}
    \begin{split}
     \frac{d}{dt}\|&\feu-\feu_h\|_{L^2(\Omega)}^2+\|\fe{F(Du)}-\F(\Du_h)\|_{L^2(\Omega)}^2\\
      &\leq c\big|(\partial_t(\feu-\feu_h),\feu-\boldsymbol{\zeta}_h)\big|
      +c\big|(\SDu-\S(\Du_h),\Du-\D\boldsymbol{\zeta}_h)\big|\\
      &\quad +c\big|(q-\mu_h,\operatorname{div}(\boldsymbol{\zeta}_h-\feu))\big|   
      +c\big|(q-\mu_h,\operatorname{div}(\feu-\feu_h))\big|\\
      &=:I_1+I_2+I_3+I_4.
    \end{split}
   \end{equation}
    $I_1$ is estimated by H\"{o}lder's inequality yielding the first
    term on the right-hand side of~\eqref{bestapp1}. For $I_2$ we use \Cref{Au4} to get
   \[I_2\leq \varepsilon\|\fe{F(Du)}-\F(\Du_h)\|_{L^2(\Omega)}^2+c_\varepsilon \|\fe{F(Du)}-\F(\D\boldsymbol{\zeta}_h)\|_{L^2(\Omega)}^2.\]
   For $I_3$ and $I_4$ we note that for a vector field $\fev$ there holds $|\operatorname{div}\fev\,|=|\operatorname{tr}(\nabla\fev)|=|\operatorname{tr}(\Dv)|\leq |\Dv|.$
   We use Young's inequality \cref{young5} for $\varphi_{|\Du|}$ to get
   \begin{equation*}
   \begin{split}
   I_3+I_4\leq& c_\varepsilon\int\limits_\Omega \big(\varphi_{|\Du|}\big)^*(|q-\mu_h|) \,dx
   +c\int\limits_\Omega \varphi_{|\Du|}(|\Du-\D\boldsymbol{\zeta}_h|) \,dx\\
   &+\varepsilon\int\limits_\Omega \varphi_{|\Du|}(|\Du-\D\feu_h|) \,dx\\
   \leq& c_\varepsilon\int\limits_\Omega \big(\varphi_{|\Du|}\big)^*(|q-\mu_h|) \,dx
   +c \|\fe{F(Du)}-\F(\D\boldsymbol{\zeta}_h)\|^2_{L^{2}(I,L^{2}(\Omega))}\\
   &+\varepsilon \|\fe{F(Du)}-\F(\D\feu_h)\|^2_{L^{2}(I,L^{2}(\Omega))}.
   \end{split}
   \end{equation*}
  Here, we used \Cref{pot15} in the last step. Choosing $\varepsilon$ sufficiently small, the assertion follows.  
  A more elaborated version of this proof for the stationary case can be found in \cite[Lemma 3.1]{BBDR}. \enddemo
    
  The terms on the right-hand side of \cref{bestapp1} can be estimated
  using the following result:
\begin{Lemma}\label{artabsch}
   Let \Cref{BP} and \Cref{io1} be fulfilled.
   Moreover, assume that $\feu_0^h\in V_h$ is given by \cref{iv} and let $\feu_h$ be the corresponding finite element solution. Then we have
  \[ \|\feu-\pihdivu\|_{L^{2}(I,L^{2}(\Omega))} + 
  \|\fe{F(Du)}-\F(\D\pihdivu)\|^2_{L^{2}(I,L^{2}(\Omega))}\leq c(K_2) h^2.\]
   For the pressure term, there holds
    \begin{equation*}
    \int\limits_I\int\limits_\Omega\big(\varphi_{|\fe{Du}|}\big)^*(|q-\Pi_h^Y q|) \,dx\, dt
      \leq c(K_2,K_3)h^{\min\{2,p'\}}.
 \end{equation*}
\end{Lemma}
\demo{}See \cite{BBDR} and \Cref{io11}.\enddemo
  Integrating \cref{bestapp1} in time produces the term
  $\|\feu_0-\feu_0^h\|_{L^2(\Omega)}$, which 
has to be estimated. Recall that 
  $\feu_0^h$ is given as the solution of \cref{iv}. 
  \begin{Lemma}\label{bestapp4}
  Let \Cref{BP} and \Cref{io1} be fulfilled.
   Moreover, assume that $\feu_0^h\in V_h$ is given by \cref{iv} and
   let $\feu_h$ be the corresponding finite element solution ensured
   by \Cref{fe11}. For $p\in[\frac{2d}{d+2},\infty)$ we have
  \begin{equation}\label{bestapp4a}
  \|\feu_0-\feu_0^h\|_{L^2(\Omega)}\leq c(K_1) h^{\min\{1,\frac{2}{p}\}}.
  \end{equation}
  \end{Lemma}
\demo{}
  Equation \cref{iv} implies the orthogonality
  \[(\fe{S(Du}_0^h)- \fe{S(Du}_0) , \D\boldsymbol{\xi}_h)=0\qquad\text{for all}\ \boldsymbol{\xi}_h\in V_h.\]
  Since $\pihdiv:W_{0,\operatorname{div}}^{1,p}(\Omega)\to V_h$, we get from \Cref{pot15}
  \begin{equation*}
  \begin{split} \|\F(\D\feu_0)-\F(\D\feu_0^h)\|_{L^2(\Omega)}^2
    &\leq c\big|(\fe{S(Du}_0)-\fe{S(Du}_0^h),\D\feu_0-\D\feu_0^h)\big|\\
    &= c\big|(\fe{S(Du}_0)-\fe{S(Du}_0^h),\D\feu_0-\D\pihdiv\feu_0)\big|.
    \end{split}
  \end{equation*}
  By \Cref{Au4}, this implies
  \[\|\F(\D\feu_0)-\F(\D\feu_0^h)\|_{L^2(\Omega)}^2\leq c\|\F(\D\feu_0)-\F(\D\pihdiv\feu_0)\|_{L^2(\Omega)}^2. \]
  Using the properties of $\pihdiv$, see \cref{io16}, we get
  \begin{equation}\label{bestapp4c}
  \|\F(\D\feu_0)-\F(\D\feu_0^h)\|_{L^2(\Omega)}^2\leq ch^2\|\nabla\F(\Du_0)\|_{L^2(\Omega)}^2\leq c(K_1)h^2.
  \end{equation}
  For the case $p\in[\frac{2d}{d+2},2)$, we may use the embedding $W_0^{1,p}(\Omega)\hookrightarrow L^2(\Omega)$, Poincar\'{e}'s and Korn's inequality, as well as \cref{Au7}, \cref{iv5} and \cref{bestapp4c} to get
  \begin{equation*}
  \begin{split}
    \|\feu_0-\feu_0^h\|_{L^2(\Omega)}^2&\leq c\|\Du_0-\Du_0^h\|_{L^p(\Omega)}^2\\
   &\leq c(K+\|\Du_0\|_{L^p(\Omega)}+\|\Du_0^h\|_{L^p(\Omega)})^{2-p} \|\F(\D\feu_0)-\F(\D\feu_0^h)\|_{L^2(\Omega)}^2\\
  &\leq c(K_1)\|\F(\D\feu_0)-\F(\D\feu_0^h)\|_{L^2(\Omega)}^2\\
  &\leq c(K_1)h^2.
  \end{split}
  \end{equation*}
  For $p\in[2,\infty),$ we get in a similar manner, this time using \cref{Au8},
  \begin{equation*}
  \|\feu_0-\feu_0^h\|_{L^2(\Omega)}^p\leq c\|\Du_0-\Du_0^h\|_{L^p(\Omega)}^p
    \leq c\|\F(\D\feu_0)-\F(\D\feu_0^h)\|_{L^2(\Omega)}^2
   \leq c(K_1)h^2,
  \end{equation*}
  which yields the assertion.
\enddemo

Now we are able to prove the main result of this section.
\begin{Theorem}\label{FEerr}
  Let \Cref{BP} and \Cref{io1} be fulfilled.
   Moreover, assume that $\feu_0^h\in V_h$ is given by \cref{iv} and
   let $\feu_h$ be the corresponding finite element solution ensured
   by \Cref{fe11}. Then, for $p\in[\frac{2d}{d+2},2)$ we have
  \begin{equation}\label{FEerra}
   \|\feu-\feu_h\|_{L^{\infty}(I,L^{2}(\Omega))}+\|\FDu-\F(\Du_h)\|_{L^{2}(I,L^{2}(\Omega))}\leq ch
  \end{equation}
   and for $p\in[2,\infty)$ we obtain
  \begin{equation}\label{FEerrb}
   \|\feu-\feu_h\|_{L^{\infty}(I,L^{2}(\Omega))}+\|\FDu-\F(\Du_h)\|_{L^{2}(I,L^{2}(\Omega))}\leq c h^{\frac{2}{p}}
  \end{equation}
  with constants depending on $K_1,\,K_2$, and $K_3$.
\end{Theorem}
\demo{}  
We choose $\boldsymbol{\zeta}_h=\pihdivu(t)$ and $\mu_h=\Pi_h^Yq(t)$ in \Cref{bestapp} 
 and integrate over $t\in I$ to get
  \begin{equation*}
   \begin{split}
   \|\feu-\feu_h&\|_{L^{\infty}(I,L^{2}(\Omega))}^2+\|\FDu-\F(\Du_h)\|_{L^{2}(I,L^{2}(\Omega))}^2 \\
     &\leq  \ c(K_1,K_2)\|\feu-\pihdivu\|_{L^{2}(I,L^{2}(\Omega))} +c  \|\fe{F(Du)}-\F(\D\pihdivu)\|^2_{L^{2}(I,L^{2}(\Omega))}\\ 
    &\ \ \ +c\int\limits_I\int\limits_\Omega\big(\varphi_{|\Du|}\big)^*(|q-\Pi_h^Yq|) \,dx \,dt +c\|\feu_0-\feu_0^h\|_{L^2(\Omega)}^2,
   \end{split}
  \end{equation*}
 where we used \Cref{BP} and \cref{rem2a} to bound $\|\partial_t(\feu-\feu_h)\|_{L^{2}(I,L^{2}(\Omega))}$. 
  Now
  \Cref{artabsch}, \Cref{bestapp4} and \Cref{BP} yield
       \begin{equation*}
   \begin{split}
   \|\feu-\feu_h&\|_{L^{\infty}(I,L^{2}(\Omega))}^2+\|\FDu-\F(\Du_h)\|_{L^{2}(I,L^{2}(\Omega))}^2 \\
     &\leq c(K_1,K_2)h^2 +c(K_2,K_3)h^{\min\{2,p'\}}+c(K_1) h^{\min\{2,\frac{4}{p}\}}.
   \end{split}
  \end{equation*}
 Since
  \[h^{\min\{2,p'\}}\leq h^{\min\{2,\frac{4}{p}\}},\]
  the theorem is proven.
\enddemo

  
  \section{The Fully Discrete Solution}\label{The fully discrete solution}
In order to numerically compute an approximate solution for the generalized Stokes equations, we still need to get rid of 
the continuity in the time variable.
Therefore, we start by dividing the time interval $I=(0,T)$ in $M$ equidistant intervals $I_n:=(t_{n-1},t_n),$ $n=1,...,M$, where $t_0=0,\ t_n=n\Delta t$, and $\Delta t=\frac{T}{M}$. For technical reasons, we assume $\Delta t<1$. The discrete time derivative is given by 
 \[d_tg^n:=\frac{g^n-g^{n-1}}{\Delta t},\qquad n=1,...,M\]
 for a sequence $(g^n)_{n=0,...,M}$ of functions $g^n\in L^1(\Omega)$.
Let $\f^n\in L^2(\Omega)$ be a suitable approximation of $\f(t_n)$ to
be specified later. The implicit scheme for the fully discrete problem
reads as follows: Given $\feU^0\in V_h$ and $\f^n\in L^2(\Omega)$ find $\feU^n\in V_h,\ n=1,...,M,$ as the solution of
    \begin{equation}\label{fd1}
    	\begin{split}
       \left(\frac{\feU^n-\feU^{n-1}}{\Delta t},\boldsymbol{\xi}_h\right)+(\fe{S(DU}^n),\D\boldsymbol{\xi}_h)=(\f^n,\boldsymbol{\xi}_h)\qquad\text{for all}\ \boldsymbol{\xi}_h\in V_h.
      \end{split}
    \end{equation}
\subsection{Existence and Regularity for the Fully Discretized Solution}

  
At first we show existence of the fully discrete solution $\feU^n$.
\begin{Theorem}\label{fd2}
  Suppose $\f^n\in L^2(\Omega),\ n=1,...,M,$ and $\feU^0\in V_h$ satisfy 
  \begin{equation}\label{fd3}
  \Delta t\sum\limits_{n=1}^M\|\f^n\|_{L^2(\Omega)}^2+\|\feU^0\|_{L^2(\Omega)}^2+\|\F(\D\feU^0)\|_{L^2(\Omega)}^2\leq K.
  \end{equation}
  Then for every $n\in\{1,...,M\}$, there exists a unique solution
  $\feU^n\in V_h$ of \cref{fd1}. If $\Delta t\leq \alpha<1$
  we obtain
  \begin{equation}\label{fd4}
  \underset{n\in\{1,...,M\}}{\sup}\|\feU^n\|_{L^2(\Omega)}^2+\Delta t\sum\limits_{n=1}^M\|\fe{F(DU}^n)\|_{L^2(\Omega)}^2\leq c(K,\alpha)
  \end{equation}
   uniformly in $M$ and $\Delta t$.
\end{Theorem}
\demo{}
The existence of $\feU^n$ follows from Brouwer's fixed point theorem.
    Setting $\boldsymbol{\xi}_h=\feU^n$ in \cref{fd1}, we get, using also the definitions of $\S$ and $\F$,  
        \[\frac{1}{2\Delta t}\|\feU^n\|_{L^2(\Omega)}^2-\frac{1}{2\Delta t}\|\feU^{n-1}\|_{L^2(\Omega)}^2+\|\fe{F(DU}^n)\|_{L^2(\Omega)}^2
    \leq|(\f^n,\feU^n)|,
    \]
 since $-(\feU^n,\feU^{n-1})=\frac{1}{2}\|\feU^n-\feU^{n-1}\|_{L^2(\Omega)}^2-\frac{1}{2}\|\feU^n\|_{L^2(\Omega)}^2-\frac{1}{2}\|\feU^{n-1}\|^2_{L^2(\Omega)}$.
    Summation from $n=1,...,l,\ l\in\{1,...,M\}$, gives
  \begin{equation*}
  \begin{split}
    \|\feU^l\|_{L^2(\Omega)}^2&+2\Delta t\sum\limits_{n=1}^l\|\fe{F(DU}^n)\|_{L^2(\Omega)}^2\\
    &\leq \|\feU^0\|_{L^2(\Omega)}^2+\Delta t\sum\limits_{n=1}^l\|\f^n\|_{L^2(\Omega)}^2
    +\Delta t\sum\limits_{n=1}^l\|\feU^n\|_{L^2(\Omega)}^2.
    \end{split}
  \end{equation*}
  Now, the discrete Gronwall lemma, see \cite{HR4},  shows
  \begin{equation*}
  \begin{split}
  \|\feU^l\|_{L^2(\Omega)}^2&+2\Delta t\sum\limits_{n=1}^l\|\fe{F(DU}^n)\|_{L^2(\Omega)}^2\\
   &\leq \exp\Big(\frac{M\Delta t}{1-\Delta t}\Big)\Big(\Delta t\sum\limits_{n=1}^M\|\f^n\|_{L^2(\Omega)}^2 +\|\feU^0\|_{L^2(\Omega)}^2\Big).
   \end{split}
  \end{equation*}
  Since $M=\frac{T}{\Delta t}$ and $\frac{1}{1-\Delta t}\leq \frac
  1{1-\alpha}$, \cref{fd4} is proven with $c(K,\alpha)=2\exp(\frac T{1-\alpha})K$.
   Uniqueness of $\feU^n$ can be proven by a similar calculation.
\enddemo
  \subsection{Error Estimates}\label{chap4err1}
Now we show error estimates between the semi-discrete solution $\feu_h(t_n)$ and the fully discrete solution $\feU^n$. To this end, we first need to establish an error equation.
For $g\in L^1(I)$, we define its mean-value on $I_n$ by
\[\overline{g}^n:=\frac{1}{\Delta t}\int\limits_{t_{n-1}}^{t_n}g(s)\,ds.\]
In \cref{fd1} we choose 
\[\f^n:=\overline{\f}^n,\qquad n=1,...,M\]
as an approximate for $\f(t_n)$. This leads to
\begin{Proposition}\label{ex}
 Let \Cref{BP} and \Cref{io1} be fulfilled.
   Moreover, assume that $\feu_0^h\in V_h$ is given by \cref{iv}. 
   Set $\f^n=\overline{\f}^n$ and $\feU^0=\feu_0^h$. Then for every $n\in\{1,...,M\},$ there exists a unique solution $\feU^n\in V_h$ of \cref{fd1} that satisfies 
  \begin{equation*}
  \underset{n\in\{1,...,M\}}{\sup}\|\feU^n\|_{L^2(\Omega)}^2+\Delta t\sum\limits_{n=1}^M\|\fe{F(DU}^n)\|_{L^2(\Omega)}^2\leq c(K_1,\alpha)
  \end{equation*}
  provided $\Delta t\leq \alpha<1$.
\end{Proposition}
\demo{}
By definition of $\overline{\f}^n$, it is clear that
\[
   \Delta t\sum\limits_{n=1}^M\|\f^n\|_{L^2(\Omega)}^2 
    =\Delta t \sum\limits_{n=1}^M\int\limits_\Omega\left|\frac{1}{\Delta t}\int\limits_{I_n}\f(s)\,ds\right|^2 \,dx
    \leq \sum\limits_{n=1}^M\int\limits_\Omega \int\limits_{I_n}|\f(s)|^2\,ds \,dx
    \leq c(K_1).
\]
\Cref{iv1} yields
  \[\|\feu_0^h\|_{L^2(\Omega)}+\|\F(\Du_0^h)\|_{L^2(\Omega)}\leq c(K_1),\]
where we used the embedding $W^{1,p}(\Omega)\hookrightarrow L^2(\Omega)$, which holds for $p\in[\frac{2d}{d+2},\infty)$.
Therefore, the requirements of \Cref{fd2} are fulfilled and the proposition is proven.
\enddemo
For the error equation, we take the mean-value of \cref{fe2} on $I_n$ and subtract \cref{fd1} to get
\begin{equation}\label{fderr1}
 (d_t(\feu_h^n-\feU^n),\boldsymbol{\xi}_h)+(\overline{\S(\Du_h)}^n-\S(\D\feU^n),\D\boldsymbol{\xi}_h)=0
\end{equation}
for all $\boldsymbol{\xi}_h\in V_h,\ n=1,...,M$. Here we use the notation
\[\feu_h^n:=\feu_h(t_n),\qquad n=1,...,M.\]
As a first step, we state a preliminary result
\begin{Lemma}\label{fderr2}
    For $n=1,...,M$ there holds
    \[\|\F(\Du_h^n)-\overline{\F(\Du_h)}^n\|_{L^2(\Omega)}^2\leq\Delta t\|\partial_t\F(\Du_h)\|_{L^2(I_n,L^2(\Omega))}^2.\]
\end{Lemma}
\demo{}
    The fundamental theorem of calculus yields
    \begin{equation*}
      \begin{split}
      \|\F(\Du_h^n)-\overline{\F(\Du_h)}^n\|_{L^2(\Omega)}^2
      &=\int\limits_\Omega\Big|\dashint\limits_{I_n}\int\limits_s^{t_n}\partial_\tau\F(\Du_h(\tau))d\tau\,ds\Big|^2 \,dx\\
      &\leq\int\limits_\Omega\Big|\int\limits_{t_{n-1}}^{t_n}\partial_\tau\F(\Du_h(\tau))d\tau\Big|^2 \,dx\\
      &\leq(t_n-t_{n-1})\|\partial_t\F(\Du_h)\|_{L^2(I_n,L^2(\Omega))}^2,
      \end{split}
    \end{equation*}
    where we used H\"{o}lder's inequality in the last step. This proves the lemma.
\enddemo
Now we are ready to prove the error estimate between the semi discretized solution $\feu_h^n=\feu_h(t_n)$ and the fully discrete solution $\feU^n$.
\begin{Proposition}\label{fderrmain}
  Let \Cref{BP} and \Cref{io1} be fulfilled.
   Moreover, assume that $\feu_0^h\in V_h$ is given by \cref{iv} and
   let $\feu_h$ be the corresponding finite element solution ensured
   by \Cref{fe11}. Set $\f^n=\overline{\f}^n$ and $\feU^0=\feu_0^h$. Then we have
  \begin{equation*}
  \underset{n\in\{1,...,M\}}{\sup}\|\feu_h^n-\feU^n\|_{L^2(\Omega)}^2+\Delta t\sum\limits_{n=1}^M\|\F(\Du_h^n)- \F(\D\feU^n)\|_{L^2(\Omega)}^2 \leq c(K_1)(\Delta t)^2
  \end{equation*}
  provided $\Delta t\leq \alpha<1$.
\end{Proposition}
\demo{}
 We choose $\boldsymbol{\xi}_h=\feu_h^n-\feU^n$ as a test function in the error equation \cref{fderr1}. After rearranging the terms and using \Cref{pot15}
  and \Cref{Au4}, we get
    \openup 1\jot
  \begin{align}
   (d_t(\feu_h^n&-\feU^n),\feu_h^n-\feU^n)+\|\F(\Du_h^n)- \F(\D\feU^n)\|_{L^2(\Omega)}^2\nonumber \\
    &\leq c\big|(\S(\Du_h^n)-\overline{\S(\Du_h)}^n,\D\feu_h^n-\D\feU^n )\big|\nonumber\\
    &=c\big|\dashint\limits_{I_n} (\S(\Du_h(t_n))-\S(\Du_h(s)),\D\feu_h(t_n)-\D\feU^n) \,ds\big|\label{fderrmain2} \\
    &\leq c_\varepsilon \dashint\limits_{I_n} \|\F(\Du_h(t_n)-\F(\Du_h(s))\|_{L^2(\Omega)}^2 \,ds
       +\varepsilon \|\F(\Du_h^n)- \F(\D\feU^n)\|_{L^2(\Omega)}^2\nonumber\\
    &\leq c_\varepsilon \|\F(\Du_h^n)-\overline{\F(\Du_h)}^n\|_{L^2(\Omega)}^2
         +\varepsilon \|\F(\Du_h^n)- \F(\D\feU^n)\|_{L^2(\Omega)}^2.\nonumber
  \end{align}
 In the first term on the left-hand side of \cref{fderrmain2}, a simple calculation yields
 \begin{equation*}
   (d_t(\feu_h^n-\feU^n),\feu_h^n-\feU^n) \geq\frac{1}{2\Delta t}\big(\|\feu_h^n-\feU^n\|_{L^2(\Omega)}^2-\|\feu_h^{n-1}-\feU^{n-1}\|_{L^2(\Omega)}^2\big).\\
  \end{equation*}
 Taking this into account and choosing $\varepsilon$ sufficiently small, we get from \cref{fderrmain2}
 \begin{equation*}
   \openup 1\jot
  \begin{split}
     \|\feu_h^n&-\feU^n\|_{L^2(\Omega)}^2-\|\feu_h^{n-1}-\feU^{n-1}\|_{L^2(\Omega)}^2 +\Delta t \|\F(\Du_h^n)- \F(\D\feU^n)\|_{L^2(\Omega)}^2\\
      &\leq c\Delta t \|\F(\Du_h^n)-\overline{\F(\Du_h)}^n\|_{L^2(\Omega)}^2
      \leq c(\Delta t)^2\|\partial_t\F(\Du_h)\|_{L^2(I_n,L^2(\Omega))}^2,
  \end{split}
 \end{equation*}
 where we also used \Cref{fderr2} in the last step. Now summation from $n=1,...,l$, taking the supremum over $l\in\{1,...,M\}$ and taking \cref{rem2a} into account yields the assertion.
 \enddemo

In order to link the continuous function $\feu$ to the fully discrete
function $\feU^n$, $n=1,...,M$, we define the piecewise-constant-in-time function
\begin{equation*}
 \hat{\feU}(t):=\begin{cases}
         &\feU^0,\qquad t=0\\
         &\feU^n,\qquad t\in I_n,\,n=1,...,M. 
        \end{cases}
\end{equation*}
Together with the results from the previous section, we get our main error estimate
\begin{Theorem}\label{geserr1}
 Let \Cref{BP} and \Cref{io1} be fulfilled.
   Moreover, assume that $\feu_0^h\in V_h$ is given by \cref{iv} and
   let $\feU^n$ be the corresponding fully discrete solution ensured
   by \Cref{ex}. For $p\in[\frac{2d}{d+2},2]$, we have
  \[\|\feu-\hat{\feU}\|_{L^{\infty}(I,L^{2}(\Omega))}+\|\FDu-\F(\D\hat{\feU})\|_{L^{2}(I,L^{2}(\Omega))}\leq c(K_1,K_2,K_3)h+c(K_1)\Delta t,\]
 and for $p\in(2,\infty)$, we have
\[\|\feu-\hat{\feU}\|_{L^{\infty}(I,L^{2}(\Omega))}+\|\FDu-\F(\D\hat{\feU})\|_{L^{2}(I,L^{2}(\Omega))}\leq c(K_1,K_2,K_3)h^{\frac{2}{p}}+c(K_1)\Delta t,\]
  provided $\Delta t\leq \alpha<1$.
\end{Theorem}
\demo{}
Define $\beta:=1$, if $p\in[\frac{2d}{d+2},2]$ and $\beta:=\frac{2}{p}$ if $p\in(2,\infty)$.
Let $\feu_h$ be the finite element solution ensured by \Cref{fe11}. Note that the fundamental theorem of calculus yields for $t\in I_n$
  \begin{equation*}
     \|\feu_h(t)-\feu_h(t_n)\|_{L^2(\Omega)}^2    
        \leq \int\limits_\Omega\Big|\int\limits_{t_{n-1}}^{t_n}\partial_\tau\feu_h(\tau,x) d\tau\Big|^2  \,dx
      \leq(\Delta t)^2\|\partial_t\feu_h\|_{L^\infty(I_n,L^2(\Omega))}^2.
\end{equation*}
This, \Cref{FEerr}, 
 and \Cref{fderrmain} yield that
  \begin{equation*}
 \begin{split}
   \underset{t\in I_n}{\sup}\|\feu-\hat{\feU}\|_{L^2(\Omega)}^2
   &\leq \underset{t\in I_n}{\sup} \|\feu-\feu_h\|_{L^2(\Omega)}^2 +\underset{t\in I_n}{\sup} \|\feu_h-\feu_h^n\|_{L^2(\Omega)}^2 +\underset{t\in I_n}{\sup} \|\feu_h^n-\feU^n\|_{L^2(\Omega)}^2\\
   & \leq c(K_1,K_2,K_3)h^{2\beta}+(\Delta t)^2\|\partial_t\feu_h\|_{L^\infty(I_n,L^2(\Omega))}^2 +c(K_1)(\Delta t)^2.
  \end{split}
  \end{equation*}
Moreover, we have
  \begin{equation*}
\begin{split}
   \|\FDu&-\F(\D\hat{\feU})\|_{L^{2}(I,L^{2}(\Omega))}^2\\
     &\leq \|\FDu-\F(\D\feu_h)\|_{L^{2}(I,L^{2}(\Omega))}^2+\sum\limits_{n=1}^M \|\F(\Du_h)-\F(\D\feu_h^n)\|_{L^2(I_n,L^2(\Omega))}^2\\
      & \quad + \sum\limits_{n=1}^M \|\F(\Du_h^n)-\F(\D\feU^n)\|_{L^2(I_n,L^2(\Omega))}^2\\
     &=:I_1+I_2+I_3.
  \end{split}
  \end{equation*}
From \Cref{FEerr}, we have 
\[I_1\leq c(K_1,K_2,K_3)h^{2\beta}\]
and \Cref{fderrmain} gives
\[   I_3= \sum\limits_{n=1}^M (t_n-t_{n-1}) \|\F(\Du_h^n)-\F(\D\feU^n)\|_{L^2(\Omega))}^2\leq c(K_1)(\Delta t)^2.
  \]
A similar argument as in \Cref{fderr2} shows
\[
   I_2=\sum\limits_{n=1}^M \int\limits_{I_n}\int\limits_\Omega|\F(\Du_h(s))-\F(\D\feu_h(t_n))|^2 \,dx \,ds
     \leq (\Delta t)^2 \|\partial_t\F(\Du_h)\|_{L^2(I,L^2(\Omega))}^2.
\]
Altogether we obtain, also using \Cref{rem2}, the assertion of the theorem.
\enddemo


%
%
%
%
%
\def\cprime{$'$}


\begin{thebibliography}{10}

\bibitem{BL1}
{\sc J.~W. Barrett and W.~B. Liu}, {\em Finite element approximation of the
  {$p$}-{L}aplacian}, Math. Comp., 61 (1993), pp.~523--537.

\bibitem{BL2}
{\sc J.~W. Barrett and W.~B. Liu}, {\em Finite element approximation of the
  parabolic {$p$}-{L}aplacian}, SIAM J. Numer. Anal., 31 (1994), pp.~413--428.

\bibitem{BL3}
{\sc J.~W. Barrett and W.~B. Liu}, {\em Quasi-norm error bounds for the finite
  element approximation of a non-{N}ewtonian flow}, Numer. Math., 68 (1994),
  pp.~437--456.

\bibitem{BBDR}
{\sc L.~Belenki, L.~C. Berselli, L.~Diening, and M.~R{\r{u}}{\v{z}}i{\v{c}}ka},
  {\em On the finite element approximation of {$p$}-{S}tokes systems}, SIAM J.
  Numer. Anal., 50 (2012), pp.~373--397.

\bibitem{BDR09}
{\sc L.~C. Berselli, L.~Diening, and M.~R{\r{u}}{\v{z}}i{\v{c}}ka}, {\em
  Optimal error estimates for a semi-implicit {E}uler scheme for incompressible
  fluids with shear dependent viscosities}, SIAM J. Numer. Anal., 47 (2009),
  pp.~2177--2202.

\bibitem{BDR}
{\sc L.~C. Berselli, L.~Diening, and M.~R{\r{u}}{\v{z}}i{\v{c}}ka}, {\em
  Existence of strong solutions for incompressible fluids with shear dependent
  viscosities}, J. Math. Fluid Mech., 12 (2010), pp.~101--132.

\bibitem{BDR14}
{\sc L.~C. Berselli, L.~Diening, and M.~R{\r{u}}{\v{z}}i{\v{c}}ka}, {\em
  Optimal error estimate for semi-implicit space-time discretization for the
  equations describing incompressible generalized {N}ewtonian fluids}, IMA J.
  Numer. Anal., 35 (2015), pp.~680--697.

\bibitem{BP}
{\sc D.~Bothe and J.~Pr{\"u}ss}, {\em {$L_P$}-theory for a class of
  non-{N}ewtonian fluids}, SIAM J. Math. Anal., 39 (2007), pp.~379--421.

\bibitem{BF}
{\sc F.~Brezzi and M.~Fortin}, {\em Mixed and hybrid finite element methods},
  vol.~15 of Springer Series in Computational Mathematics, Springer-Verlag, New
  York, 1991.

\bibitem{DER}
{\sc L.~Diening, C.~Ebmeyer, and M.~R{\r{u}}{\v{z}}i{\v{c}}ka}, {\em Optimal
  convergence for the implicit space-time discretization of parabolic systems
  with {$p$}-structure}, SIAM J. Numer. Anal., 45 (2007), pp.~457--472
  (electronic).

\bibitem{DNR}
{\sc L.~Diening, P.~N{\"a}gele, and M.~R{\r{u}}{\v{z}}i{\v{c}}ka}, {\em
  Monotone operator theory for unsteady problems in variable exponent spaces},
  Complex Var. Elliptic Equ., 57 (2012), pp.~1209--1231.

\bibitem{DPR}
{\sc L.~Diening, A.~Prohl, and M.~R{\r{u}}{\v{z}}i{\v{c}}ka}, {\em On
  time-discretizations for generalized {N}ewtonian fluids}, in Nonlinear
  problems in mathematical physics and related topics, {II}, vol.~2 of Int.
  Math. Ser. (N. Y.), Kluwer/Plenum, New York, 2002, pp.~89--118.

\bibitem{DPR06}
{\sc L.~Diening, A.~Prohl, and M.~R{\r{u}}{\v{z}}i{\v{c}}ka}, {\em
  Semi-implicit {E}uler scheme for generalized {N}ewtonian fluids}, SIAM J.
  Numer. Anal., 44 (2006), pp.~1172--1190 (electronic).

\bibitem{DR3}
{\sc L.~Diening and M.~R{\r{u}}{\v{z}}i{\v{c}}ka}, {\em Strong solutions for
  generalized {N}ewtonian fluids}, J. Math. Fluid Mech., 7 (2005),
  pp.~413--450.

\bibitem{DR2}
{\sc L.~Diening and M.~R{\r{u}}{\v{z}}i{\v{c}}ka}, {\em Interpolation operators
  in {O}rlicz-{S}obolev spaces}, Numer. Math., 107 (2007), pp.~107--129.

\bibitem{Sarah}
{\sc S.~Eckstein}, {\em On the full space-time discretization of the
  generalized Stokes systems: The Dirichlet case}, PhD thesis,
  Albert-Ludwigs-Universit{\"a}t Freiburg im Breisgau, 2016.

\bibitem{GaGr}
{\sc H.~Gajewski, K.~Gr{\"o}ger, and K.~Zacharias}, {\em Nichtlineare
  {O}peratorgleichungen und {O}peratordifferentialgleichungen},
  Akademie-Verlag, Berlin, 1974.
\newblock Mathematische Lehrb{\"u}cher und Monographien, II. Abteilung,
  Mathematische Monographien, Band 38.

\bibitem{GL}
{\sc V.~Girault and J.-L. Lions}, {\em Two-grid finite-element schemes for the
  steady {N}avier-{S}tokes problem in polyhedra}, Port. Math. (N.S.), 58
  (2001), pp.~25--57.

\bibitem{GR}
{\sc V.~Girault and P.-A. Raviart}, {\em Finite element methods for
  {N}avier-{S}tokes equations}, vol.~5 of Springer Series in Computational
  Mathematics, Springer-Verlag, Berlin, 1986.
\newblock Theory and algorithms.

\bibitem{GS}
{\sc V.~Girault and L.~R. Scott}, {\em A quasi-local interpolation operator
  preserving the discrete divergence}, Calcolo, 40 (2003), pp.~1--19.

\bibitem{HR1}
{\sc J.~G. Heywood and R.~Rannacher}, {\em Finite element approximation of the
  nonstationary {N}avier-{S}tokes problem. {I}. {R}egularity of solutions and
  second-order error estimates for spatial discretization}, SIAM J. Numer.
  Anal., 19 (1982), pp.~275--311.

\bibitem{HR2}
{\sc J.~G. Heywood and R.~Rannacher}, {\em Finite element approximation of the
  nonstationary {N}avier-{S}tokes problem. {II}. {S}tability of solutions and
  error estimates uniform in time}, SIAM J. Numer. Anal., 23 (1986),
  pp.~750--777.

\bibitem{HR3}
{\sc J.~G. Heywood and R.~Rannacher}, {\em Finite element approximation of the
  nonstationary {N}avier-{S}tokes problem. {III}. {S}moothing property and
  higher order error estimates for spatial discretization}, SIAM J. Numer.
  Anal., 25 (1988), pp.~489--512.

\bibitem{HR4}
{\sc J.~G. Heywood and R.~Rannacher}, {\em Finite-element approximation of the
  nonstationary {N}avier-{S}tokes problem. {IV}. {E}rror analysis for
  second-order time discretization}, SIAM J. Numer. Anal., 27 (1990),
  pp.~353--384.

\bibitem{Krasno}
{\sc M.~A. Krasnosel{\cprime}ski{\v\i} and J.~B. Ruticki{\v\i}}, {\em Convex
  functions and {O}rlicz spaces}, P. Noordhoff Ltd., Groningen, 1961.

\bibitem{La}
{\sc O.~A. Ladyzhenskaya}, {\em The mathematical theory of viscous
  incompressible flow}, Second English edition, revised and enlarged.
  Mathematics and its Applications, Vol. 2, Gordon and Breach, Science
  Publishers, New York-London-Paris, 1969.

\bibitem{Li}
{\sc J.-L. Lions}, {\em Quelques m\'ethodes de r\'esolution des probl\`emes aux
  limites non lin\'eaires}, Dunod; Gauthier-Villars, Paris, 1969.

\bibitem{MNRR}
{\sc J.~M{\'a}lek, J.~Ne{\v{c}}as, M.~Rokyta, and
  M.~R{\r{u}}{\v{z}}i{\v{c}}ka}, {\em Weak and measure-valued solutions to
  evolutionary {PDE}s}, vol.~13 of Applied Mathematics and Mathematical
  Computation, Chapman \& Hall, London, 1996.

\bibitem{MR}
{\sc J.~M{\'a}lek and K.~R. Rajagopal}, {\em Mathematical issues concerning the
  {N}avier-{S}tokes equations and some of its generalizations}, in Evolutionary
  equations. {V}ol. {II}, Handb. Differ. Equ., Elsevier/North-Holland,
  Amsterdam, 2005, pp.~371--459.

\bibitem{PR}
{\sc A.~Prohl and M.~R{\r{u}}{\v{z}}i{\v{c}}ka}, {\em On fully implicit
  space-time discretization for motions of incompressible fluids with
  shear-dependent viscosities: the case {$p\leq 2$}}, SIAM J. Numer. Anal., 39
  (2001), pp.~214--249 (electronic).

\bibitem{DR1}
{\sc M.~R{\r{u}}{\v{z}}i{\v{c}}ka and L.~Diening}, {\em Non-{N}ewtonian fluids
  and function spaces}, in N{AFSA} 8---{N}onlinear analysis, function spaces
  and applications. {V}ol. 8, Czech. Acad. Sci., Prague, 2007, pp.~94--143.

\bibitem{SZ}
{\sc L.~R. Scott and S.~Zhang}, {\em Finite element interpolation of nonsmooth
  functions satisfying boundary conditions}, Math. Comp., 54 (1990),
  pp.~483--493.

\bibitem{Ze2}
{\sc E.~Zeidler}, {\em Nonlinear functional analysis and its applications.
  {II}/{B}}, Springer-Verlag, New York, 1990.
\newblock Nonlinear monotone operators.

\end{thebibliography}
\end{document}